\documentclass{amsart}

\newcommand{\ab}{\textnormal{ab}}

\newcommand{\diag}{\textnormal{diag}}

\newcommand{\lcm}{\textnormal{lcm}}

\newcommand{\Char}{\textnormal{Char}}
\newcommand{\Cl}{\textnormal{Cl}}

\newcommand{\Exp}{\textnormal{Exp}}
\newcommand{\Gal}{\textnormal{Gal}}
\newcommand{\GL}{\textnormal{GL}}

\newcommand{\Inf}{\textnormal{Inf}}
\newcommand{\Irr}{\textnormal{Irr}}
\newcommand{\M}{\textnormal{M}}

\newcommand{\PSp}{\textnormal{PSp}}
\newcommand{\PSL}{\textnormal{PSL}}
\newcommand{\PSU}{\textnormal{PSU}}

\renewcommand{\S}{\textnormal{S}}
\newcommand{\Sp}{\textnormal{Sp}}

\newcommand{\SL}{\textnormal{SL}}

\newcommand{\SU}{\textnormal{SU}}

\newcommand{\Z}{\textnormal{Z}}

\usepackage{algorithm}
\usepackage{algpseudocode}
\usepackage{amsmath}
\usepackage{amssymb}
\usepackage{amsthm}
\usepackage[backend=biber, style=numeric, abbreviate=true, giveninits=true]{biblatex}
\usepackage{color}
\usepackage{diagbox}
\usepackage[a4paper, left=2cm, right=2cm, top=2cm, bottom=2cm]{geometry}
\usepackage{graphicx}
\usepackage{hyperref}
\usepackage{listings}
\usepackage{mathrsfs}
\usepackage{multirow}
\usepackage{paralist}
\usepackage{wasysym}

\newtheorem{theorem}{Theorem}[section]
\newtheorem{lemma}[theorem]{Lemma}
\newtheorem{property}[theorem]{Property}
\newtheorem{proposition}[theorem]{Proposition}
\newtheorem{corollary}[theorem]{Corollary}
\newtheorem{conjecture}[theorem]{Conjecture}

\theoremstyle{definition}
\newtheorem{definition}[theorem]{Definition}
\newtheorem{example}[theorem]{Example}

\newtheorem{remark}[theorem]{Remark}

\numberwithin{equation}{section}

\AtEveryBibitem{
\clearfield{doi}
\clearfield{isbn}
\clearfield{issn}
}

\DeclareFieldFormat[article]{title}{#1}
\DeclareFieldFormat[article]{url}{}
\DeclareFieldFormat[book]{url}{}
\DeclareFieldFormat[incollection]{title}{#1}
\DeclareFieldFormat[thesis]{title}{#1}

\begin{document}

\title{\textbf{On a Conjecture of Navarro and Tiep on Character Fields}}

\author{\Large{Marco Albert}}
\address{FB Mathematik, RPTU Kaiserslautern, Postfach 3049, 67653 Kaiserslautern, Germany.}
\email{albertm@rptu.de}

\begin{abstract}
In 2021, Navarro and Tiep proposed a conjecture on character fields of finite quasi-simple groups. We develop some theory on sums of roots of unity and use this theory to prove the conjecture for some infinite families of finite quasi-simple groups with known character table. We then use the classification of the irreducible complex characters of the finite general linear groups developed by Green to obtain some partial results about the conjecture for the finite general and special linear groups in arbitrary dimension.
\end{abstract}

\maketitle

\section{Introduction}

Let $G$ be a finite group and let $\chi\in\Char(G)$, the set of complex characters of $G$.
The \textbf{field of values} $\mathbb{Q}(\chi)$ of $\chi$ is defined as the smallest field containing all values of $\chi$, i.e., the smallest extension field of $\mathbb{Q}$ such that $\chi(g)\in\mathbb{Q}(\chi)$ for all $g\in G$.
Now the value $\chi(g)$ is a sum of $o(g)$-th roots of unity, so $\mathbb{Q}(\chi)$ is a subfield of the \textbf{$n$-th cyclotomic field} $\mathbb{Q}_n$ for $n=|G|$ which is defined as the smallest extension field of $\mathbb{Q}$ containing a primitive $n$-th root of unity.
The degree $[\mathbb{Q}_n:\mathbb{Q}]$ of the $n$-th cyclotomic field over $\mathbb{Q}$ is equal to $\varphi(n)$, where $\varphi$ is Euler's totient function by \cite[Theorem 2.5]{Was}.

Apart from the degree $\chi(1)$ of a character $\chi$ there is another important number associated to $\chi$ which is its \textbf{conductor}.
For an abelian field extension $\mathbb{F}/\mathbb{Q}$, i.e., a field extension with abelian Galois group, the conductor $f_\mathbb{F}\in\mathbb{Z}_{\geq1}$ of $\mathbb{F}$ is defined as the smallest positive integer such that $\mathbb{F}\subseteq\mathbb{Q}_{f_\mathbb{F}}$.
Similarly, for some $c\in\mathbb{C}$ the conductor of $c$ is defined as $f_c:=f_{\mathbb{Q}(c)}$ if $\mathbb{Q}(c)/\mathbb{Q}$ is abelian and for a finite group $G$ and some $\chi\in\Char(G)$ the conductor of $\chi$ is defined as $f_\chi:=f_{\mathbb{Q}(\chi)}$.

It is of great interest to study the fields of values of the irreducible complex characters of $G$.
In particular, given any field extension $\mathbb{F}/\mathbb{Q}$ one might ask, whether or not there is a finite group $G$ and some $\chi\in\Irr(G)$, the set of irreducible complex chracters of $G$, such that $\mathbb{Q}(\chi)=\mathbb{F}$.
This problem has been solved in \cite[Theorem 2.2]{NaTi} namely such a group $G$ and $\chi\in\Irr(G)$ exist if and only if $\mathbb{F}/\mathbb{Q}$ is an abelian field extension.

However, it gets more difficult if we take some prime $p$ and assume that $p$ does not divide the degree of $\chi$.
For a finite group $G$ and some prime $p$ we write $\Irr_p(G)$ for the set of irreducible complex characters of $G$ of degree divisible by $p$ and $\Irr_{p'}(G)$ for the set of irreducible complex characters of degree not divisible by $p$.
So in this case the goal is to determine the set
\[\mathcal{F}_p:=\{\mathbb{Q}(\chi)\ |\ \chi\in\Irr_{p'}(G),\ G\text{ a finite group}\}.\]
It has been shown in \cite{NaTi} that given any abelian extension field $\mathbb{F}$ of $\mathbb{Q}$ with conductor $2^am$ for some $a\in\mathbb{Z}_{\geq0}$ and odd $m$ there is a finite group $G$ and some $\chi\in\Irr_{2'}(G)$ such that $\mathbb{Q}(\chi)=\mathbb{F}$ if and only if $\mathbb{Q}_{2^a}\subseteq\mathbb{F}$ which solves the problem for the case $p=2$.
To prove this statement, one needs to show that for any such field $\mathbb{F}$ there is a finite group $G$ and $\chi\in\Irr_{2'}(G)$ with $\mathbb{Q}(\chi)=\mathbb{F}$ and that no finite group $G$ and $\chi\in\Irr_{2'}(G)$ can be chosen with $\mathbb{Q}_{2^a}\nsubseteq\mathbb{F}$.
This has been done in \cite{NaTi} in the following two theorems.

\begin{theorem}[{\cite[Theorem A1]{NaTi}}]
Let $\chi\in\Irr_{2'}(G)$ have conductor $2^am$ for some $a\in\mathbb{Z}_{\geq0}$ and odd $m$.
Then $\mathbb{Q}_{2^a}\subseteq\mathbb{Q}(\chi)$.
\end{theorem}

\begin{theorem}[{\cite[Theorem A2]{NaTi}}]
Let $\mathbb{F}/\mathbb{Q}$ be an abelian field extension with conductor $2^am$ for some $a\in\mathbb{Z}_{\geq0}$ and odd $m$.
If $\mathbb{Q}_{2^a}\subseteq\mathbb{F}$, then there exist a finite group $G$ and $\chi\in\Irr_{2'}(G)$ such that $\mathbb{Q}(\chi)=\mathbb{F}$.
\end{theorem}

On the other hand, if $p$ is an odd prime, then the set $\mathcal{F}_p$ has been characterised in \cite{NaTi} in a similar way, but only modulo the following conjecture on the values of irreducible complex characters of finite quasi-simple groups, i.e., perfect groups $G$ such that $G/\Z(G)$ is simple.

\begin{property}
Let $p$ be an odd prime, let $G$ be a finite group and let $\chi\in\Irr_{p'}(G)$ with conductor $p^am$ for some $a\in\mathbb{Z}_{\geq0}$ and $m\in\mathbb{Z}_{\geq1}$ such that $p$ does not divide $m$.
We say that Property 1.3 holds for $G,\ \chi$ and $p$ if there exists a $p$-element $g\in G$ such that $p$ does not divide $[\mathbb{Q}_{p^a}:\mathbb{Q}(\chi(g))]$.
Further, we say that $G$ has Property 1.3 if it holds for $G,\ \chi$ and $p$ for every choice of an odd prime $p$ and every $\chi\in\Irr_{p'}(G)$.
\end{property}

\begin{conjecture}[{\cite[Conjecture B3]{NaTi}}]
All finite quasi-simple groups have Property 1.3.
\end{conjecture}

So if Conjecture 1.4 holds, then together with Theorem 1.1 and Theorem 1.2 the set $\mathcal{F}_p$ is characterised for an arbitrary prime $p$ by the following two theorems.

\begin{theorem}[{\cite[Theorem B1]{NaTi}}]
Suppose that Conjecture 1.4 holds.
Let $p$ be a prime, let $G$ be a finite group and let $\chi\in\Irr_{p'}(G)$ with conductor $p^am$ such that $p$ does not divide $m$ and $a\in\mathbb{Z}_{\geq0}$.
Then $[\mathbb{Q}_{p^a}:(\mathbb{Q}(\chi)\cap\mathbb{Q}_{p^a})]$ is not divisible by $p$.
\end{theorem}

\begin{theorem}[{\cite[Theorem B2]{NaTi}}]
Let $\mathbb{F}/\mathbb{Q}$ be an abelian field extension with conductor $p^am$ such that $p$ does not divide $m$ and $a\in\mathbb{Z}_{\geq0}$.
If $p$ does not divide $[\mathbb{Q}_{p^a}:(\mathbb{F}\cap\mathbb{Q}_{p^a})]$, then there exist a finite group $G$ and $\chi\in\Irr_{p'}(G)$ such that $\mathbb{Q}(\chi)=\mathbb{F}$.
\end{theorem}

Note that Theorem 1.6 does not require that Conjecture 1.4 holds, so if $\mathbb{F}/\mathbb{Q}$ is any abelian field extension with $f_\mathbb{F}=p^am$ for some prime $p$, some $a\in\mathbb{Z}_{\geq0}$ and some $m\in\mathbb{Z}_{\geq1}$ not divisible by $p$ such that $p\nmid[\mathbb{Q}_{p^a}:(\mathbb{F}\cap\mathbb{Q}_{p^a})]$, then $\mathbb{F}\in\mathcal{F}_p$.
It is just still unknown whether or not $\mathcal{F}_p$ also contains fields $\mathbb{F}$ such that $p$ divides $[\mathbb{Q}_{p^a}:(\mathbb{F}\cap\mathbb{Q}_{p^a})]$.

The main concern of this paper is Property 1.3.
We develop some results on it and prove that it holds for some infinite families of finite groups.
The results of this paper and its structure are as follows.

In Section 2 we develop some theory which helps to investigate whether or not a given group $G$ has Property 1.3.
There we prove the first main theorem of this paper, Theorem 2.5, which is a result on sums of complex $p^k$-th roots of unity for some odd prime $p$ and some $k\in\mathbb{Z}_{\geq2}$.
We then use Theorem 2.5 to give a statement which seems to be weaker than Property 1.3, but turns out to be equivalent to it.

In Section 3 we then use Theorem 2.5 to prove that some infinite series of finite quasi-simple groups with known character table have Property 1.3.

In Section 4 we use the parameterisation of the irreducible characters of the finite general linear groups developed by Green in \cite{Gre} to prove Property 1.3 for the finite general linear groups under some further assumptions (see Theorem 1.7).
Then we use Clifford Theory to prove the second main theorem of this paper, Theorem 4.6, which describes under which assumptions on a group $G$ a normal subgroup $N\unlhd G$ has Property 1.3.
Finally, we use Theorem 4.6 and the results on the finite general linear groups to prove that the finite special linear groups have Property 1.3 under similar assumptions as for the general linear groups.

Summarising all of our results on Property 1.3, we then get the third main theorem of this paper which is the following.

\begin{theorem}
Let $G$ be a finite group, let $q$ be a prime power and let $p$ be an odd prime.
Further, let $\chi\in\Irr(G)$ and set $S:=G/\Z(G)$.
Then Property 1.3 holds for $G,\ \chi$ and $p$ in all of the following cases.
\begin{itemize}
\item $S\in\{\PSL_2(q),\PSL_3(q),\PSU_3(q),\PSp_4(q),G_2(q),{^3}D_4(q)\}$ and $q$ is any prime power
\item $S\in\{\PSp_6(2^n),{^2}B_2(2^{2n+1}),{^2}G_2(3^{2n+1})\}$ for some $n\in\mathbb{Z}_{\geq1}$
\item $S={^2}F_4(2)'$
\item $S={^2}F_4(q)$, $q=2^{2n+1}$ for some $n\in\mathbb{Z}_{\geq1}$ and $\chi$ does not lie in one of the ten character families not given in \cite{CHEVIE}
\item $G=\GL_n(q),\ n\in\mathbb{Z}_{\geq1},\ p>n$ and $p|(q-1)$
\item $S=\PSL_n(q),\ n\in\mathbb{Z}_{\geq1},\ p>n(n-1)$ and $p|(q-1)$
\end{itemize}
\end{theorem}

\noindent\textbf{Acknowledgement:} This paper is a modified version of the author's master's thesis under the supervision of Prof. Dr. Gunter Malle.

\section{General theory}

Throughout let $p$ be an odd prime and let $q$ be a power of a different prime than $p$.
Moreover, let $G$ be a finite quasi-simple group and set $S:=G/\Z(G)$.
If $n=m\cdot p^a\in\mathbb{Z}_{\geq1}$ for some $a\in\mathbb{Z}_{\geq0}$ and $m\in\mathbb{Z}_{\geq1}$ such that $p\nmid m$, we set $n_p:=p^a$ and $n_{p'}:=m$.
Further, for $g\in G$ we set $g_p:=g^{o(g)_{p'}}$ and $g_{p'}:=g^{o(g)_p}$.

In this section we prove some basic tools which we need to prove Property 1.3 for the mentioned groups.
First we start with an easy case.

\begin{lemma}
Property 1.3 holds if $a\in\{0,1\}$.
In particular it holds for $G$ and $p$ if $p^2\nmid\Exp(G)$.
\begin{proof}
Let $\chi\in\Irr_{p'}(G)$.
We choose the $1$-element $\iota$ of $G$ as the $p$-element and get
\[[\mathbb{Q}_{p^a}:\mathbb{Q}(\chi(\iota))]=[\mathbb{Q}:\mathbb{Q}]=1\]
if $a=0$ respectively
\[[\mathbb{Q}_{p^a}:\mathbb{Q}(\chi(\iota))]=[\mathbb{Q}_{p}:\mathbb{Q}]=\varphi(p)=p-1,\]
where $\varphi$ is Euler's totient function, if $a=1$ both of which are not divisible by $p$.
\end{proof}
\end{lemma}

Since the character value $\chi(g)$ of any character $\chi$ of $G$ on any $g\in G$ is a sum of $o(g)$-th roots of unity, we need to study some properties of sums of roots of unity.
In particular, it is important to know when the conductors of $\chi(g)$ and $\chi(g_p)$ have the same $p$-part.
This does not always happen as shown by the example $G=\SL_6(2)$ for $p=3$.
By \cite[text below Conjecture B3]{NaTi} this is an example for a group $G$ and a prime $p$ which shows that the conclusion of Property 1.3 does not hold if $p|\chi(1)$, so there is some $g\in G$ such that $f_{\chi(g)}$ and $f_{\chi}$ have the same $p$-part, but the $p$-part of $f_{\chi(g_p)}$ is strictly smaller.
However, Theorem 2.5 shows this property for the cases we want to consider in this paper.
We start with some definitions on sums of roots of unity.

\begin{definition}
Let $n\in\mathbb{Z}_{\geq1}$ and let $z\in\mathbb{C}$ be a primitive $n$-th root of unity.
We define a \textbf{sum of $n$-th roots of unity} of length $l$ to be an $l$-tuple $s=(z^{a_1},\dots,z^{a_l})$ for some $a_1,\dots,a_l\in\mathbb{Z}$ and we identify two such tuples as the same if their entries coincide up to permutation.
We define its \textbf{value} to be
\[|s|:=z^{a_1}+\dots+z^{a_l}\]
and we call $s$ \textbf{vanishing sum} if $|s|=0$.
Further, we call the subsequences of $s$ \textbf{subsums} of $s$ and we call $s$ \textbf{minimal} if it has no non-empty proper subsum which is vanishing, i.e., if for every $\emptyset\ne I\subset\{1,\dots,l\}$ we have $\sum\limits_{i\in I}z^{a_i}\ne0$.
\end{definition}

This way of writing sums of roots of unity allows us to distinguish between two sums of roots of unity if their values are the same but their entries are not.
We need this in Theorem 2.5, but first we give an important result about minimal vanishing sums of roots of unity from \cite{LaLe}.

\begin{proposition}[{\cite[Corollary 3.4]{LaLe}}]
Let $n=p_1^ap_2^b$, where $a,b\in\mathbb{Z}_{\geq0}$ and $p_1,p_2$ are primes.
Then, up to a rotation, i.e., up to multiplication with an $n$-th root of unity, the only non-empty minimal vanishing sums of $n$-th roots of unity are $(1,z,\dots,z^{p_1-1})$ and $(1,y,\dots,y^{p_2-1})$ for a primitive $p_1$-th root of unity $z\in\mathbb{C}$ and a primitive $p_2$-th root of unity $y\in\mathbb{C}$.
\end{proposition}

One might ask if the conclusion of Proposition 2.3 also holds if $n$ is divisible by more than two primes.
However, the following example shows that this is not the case.

\begin{example}
Let $z\in\mathbb{C}$ be a root of unity of order $30=2\cdot 3\cdot 5$.
Set
\[s=(z^5,z^6,z^{12},z^{18},z^{24},z^{25}).\]
Then $s$ is a vanishing sum since
\[|s|=(z^6+z^{12}+z^{18}+z^{24})+(z^5+z^{25})=-1+1=0,\]
but $s$ has no non-empty proper vanishing subsum, so $s$ is minimal.
\end{example}

Let $s=(z^{a_1},\dots,z^{a_l})$ be a sum of roots of unity as in Definition 2.2 and let $m\in\mathbb{Z}_{\geq1}$.
We set
\begin{align*}
\lambda(s,m)&:=\text{number of entries of }s\text{ of order }m\\
\Tilde{\lambda}(s,m)&:=\text{number of different entries of }s\text{ of order }m
\end{align*}
and we write $s^{(m)}$ for the subsum of $s$ consisting of all entries of $s$ of order $m$.
Further, to simplify notation, we set $\mathbb{Q}(s):=\mathbb{Q}(|s|)$ and for the conductor of $|s|$ we set $f_s:=f_{|s|}$.

Now recall that $p$ is an odd prime.

\begin{theorem}
Let $k\in\mathbb{Z}_{\geq2}$ and let $z\in\mathbb{C}$ be a primitive $p^k$-th root of unity.
Further, let $s:=(z^{a_1},\dots,z^{a_l})$ be a sum of $p^k$-th roots of unity of length $l$ with $a_1,\dots,a_l\in\mathbb{Z}$.
\begin{compactenum}[(a)]
\item The relation $s_1\sim s_2:\Leftrightarrow|s_1|=|s_2|$ on the set of sums of $p^k$-th roots of unity is an equivalence relation.
Further, each equivalence class has a unique minimal representative and if $|s_1|=|s_2|$, then also $|s_1^{(p^k)}|=|s_2^{(p^k)}|$.
\item We have
\[f_s=p^k\Longleftrightarrow|s^{(p^k)}|\ne0.\]
In particular if $p\nmid\lambda(s,p^k)$ or $p\nmid\Tilde{\lambda}(s,p^k)$, then $f_s=p^k$.\\
\item We have
\[p|[\mathbb{Q}_{p^k}:\mathbb{Q}(s)]\Longleftrightarrow|s^{(p^k)}|=0.\]
In particular if $p\nmid\lambda(s,p^k)$ or $p\nmid\Tilde{\lambda}(s,p^k)$, then $p\nmid[\mathbb{Q}_{p^k}:\mathbb{Q}(s)]$.
\end{compactenum}
\begin{proof}
\begin{compactenum}[(a)]
\item It is clear that $"\sim"$ is an equivalence relation.
Now let $s_1$ and $s_2$ be sums of $p^k$-th roots of unity with $|s_1|=|s_2|$.
To prove the other two claims, it is enough to show that\\
(1) the equivalence class of $s_1$ has a minimal representative $s_0$ and $|s_0^{(p^k)}|=|s_1^{(p^k)}|$,\\
(2) the minimal representative of the class of $s_1$ is unique.

We start with (1).
By Proposition 2.3 any non-empty minimal vanishing subsum $t_1$ of $s_1$ is of the form
\[t_1=(z^a,z^{a+p^{k-1}},\dots,z^{a+p^{k-1}(p-1)})\]
for some $a\in\mathbb{Z}$.
Now if $p|a$, then the exponent of every entry of $t_1$ is divisible by $p$ since $k\geq2$, so none of these entries are primitive.
On the other hand, if $p\nmid a$, then the exponent of every entry of $t_1$ is not divisible by $p$ since $k\geq2$, so all of these entries are primitive.
In both cases we get a subsum $s'_1$ of $s_1$ by taking the entries of $s_1$ not contained in $t_1$ and since either all of the entries of $t_1$ are also entries of ${s'_1}^{(p^k)}$ or none of them are, we get $|{s'_1}^{(p^k)}|=|s_1^{(p^k)}|$.
We repeat this procedure until $s'_1$ is minimal.

For the proof of (2) assume that $s_1$ and $s_2$ are minimal.
Let $y\in\mathbb{C}$ be a primitive $2p^k$-th root of unity such that $y^2=z$ and write $s_1=(z^{b_1},\dots,z^{b_{l_1}})$ and $s_2=(z^{c_1},\dots,z^{c_{l_2}})$ with $b_1,\dots,b_{l_1},c_1,\dots,c_{l_2}\in\{0,\dots,p^k-1\}$.
Then $-y^a=y^{a+p^k}$ for all $a\in\mathbb{Z}$ and we get
\[y^{2b_1}+\dots+y^{2b_{l_1}}+y^{2c_1+p^k}+\dots+y^{2c_{l_2}+p^k}=0,\]
so $t:=(y^{2b_1},\dots,y^{2b_{l_1}},y^{2c_1+p^k},\dots,y^{2c_{l_2}+p^k})$ is a vanishing sum.
Applying Proposition 2.3 we get that, up to a rotation, the non-empty minimal vanishing subsums of $t$ are $(1,y^{p^k})$ and $(1,y^{2p^{k-1}},\dots,y^{2p^{k-1}(p-1)})$.
Now let $t'$ be a non-empty minimal vanishing subsum of $t$, so the length of $t'$ is either $2$ or $p$.
If $t'$ has length $p$, then either all the exponents of entries of $t'$ are even or all of them are odd.
The even ones are precisely the entries of $s_1=(y^{2b_1},\dots,y^{2b_{l_1}})$ and the odd ones are the entries of $(y^{2c_1+p^k},\dots,y^{2c_{l_2}+p^k})$.
Now by assumption $s_1$ and $s_2$ are minimal which shows that $t'$ cannot have length $p$, so every non-empty minimal vanishing subsum of $t$ has length $2$ and therefore it consists of $(l_1+l_2)/2$ minimal vanishing subsums of length $2$.
Such a subsum $t'$ consists of two entries, one with even exponent and one with odd exponent, so it has the form $t'=(y^{2b_i},y^{2c_j+p^k})$ for some $i\in\{1,\dots,l_1\},\ j\in\{1,\dots,l_2\}$ and $l_1=l_2$ follows.
This leads to the equation
\[0=y^{2b_i}+y^{2c_j+p^k}=y^{2b_i}-y^{2c_j}=z^{b_i}-z^{c_j},\]
so $b_i=c_j$.
It follows that $(c_1,\dots,c_{l_2})$ is a permutation of $(b_1,\dots,b_{l_1})$ and then $s_1=s_2$ which proves the claim.
\item First assume that $|s^{(p^k)}|=0$.
In this case $|s|$ is the value of a sum of $p^{k-1}$-th roots of unity, so $f_s<p^k$.

For the other direction assume that $|s^{(p^k)}|\ne0$.
By part (a) we can assume that $s$ is minimal since $f_s$ only depends on $|s|$ and without loss of generality we can assume that $z^{a_1}$ is a primitive $p^k$-th root of unity since $\lambda(s,p^k)\geq1$ by assumption.
Clearly $f_s\leq p^k$ since $s$ is a sum of $p^k$-th roots of unity.
Assume that $f_s<p^k$.
Then \cite[Theorem 3.4]{Con} implies that 
$\mathbb{Q}_{f_s}\subseteq\mathbb{Q}_{p^{k-1}}$, so $|s|\in\mathbb{Q}_{p^{k-1}}$ and since the set $\{1,z^p,\dots,z^{p(p^{k-1}-1)}\}$ generates $\mathbb{Q}_{p^{k-1}}$ as a $\mathbb{Q}$-vector space, we can write
\[|s|=\frac{d_0}{e_0}+\frac{d_1}{e_1}z^p+\dots+\frac{d_{p^{k-1}-1}}{e_{p^{k-1}-1}}z^{p(p^{k-1}-1)}\]
for some $d_0,\dots,d_{p^{k-1}-1}\in\mathbb{Z}$ and some $e_0,\dots,e_{p^{k-1}-1}\in\mathbb{Z}_{\geq1}$.
Set $e:=\lcm(e_0,\dots,e_{p^{k-1}-1})$.
Then by multiplying with $e$ we get
\[e|s|=d_0'+d_1'z^p+\dots+d_{p^{k-1}-1}'z^{p(p^{k-1}-1)}\]
for some $d_0',\dots,d_{p^{k-1}-1}'\in\mathbb{Z}$.
So $e|s|$ is of the form
\[e|s|=z^{b_1}+\dots+z^{b_m}-z^{c_1}-\dots-z^{c_r}\] for some 
$m,r\in\mathbb{Z}$ and some $b_1,\dots,b_m,c_1,\dots,c_r\in\mathbb{Z}$ such that $p|b_i$ and $p|c_j$ for all $i\in\{1,\dots,m\}$ and $j\in\{1,\dots,r\}$.
Now recall that $p$ is an odd prime and let $y\in\mathbb{C}$ be a primitive $2p^k$-th root of unity such that $y^2=z$.
Then again $-y^a=y^{a+p^k}$ for all $a\in\mathbb{Z}$ and by using $|s|=z^{a_1}+\dots+z^{a_l}$ we get
\[(*)\ e(y^{2a_1}+\dots+y^{2a_l})+y^{2c_1}+\dots+y^{2c_r}+y^{2b_1+p^k}+\dots+y^{2b_m+p^k}=0\]
which is the value of a vanishing sum of length $el+m+r$.
Since $2p^k$ is the product of two prime powers, we can apply Proposition 2.3
to get that, up to a rotation, the non-empty minimal vanishing subsums of the left hand side of $(*)$ are $(1,y^{p^k})$ and $(1,y^{2p^{k-1}},\dots,y^{2p^{k-1}(p-1)})$.
Now $y^{2a_1}$ was assumed to be a primitive $p^k$-th root of unity and it has to be an entry of a minimal vanishing subsum $s'$ of the left hand side of $(*)$.
So either $s'=(y^{2a_1},y^{2a_1+p^k})$ or $s'=(y^{2a_1},y^{2a_1+2p^{k-1}},\dots,y^{2a_1+2p^{k-1}(p-1)})$.
Since $k\geq2$, in both cases $p$ does not divide the exponent of any entry of $s'$ and since its entries are pairwise distinct, it must be a subsum of $(y^{2a_1},\dots,y^{2a_l})=s$ which is a contradiction to $s$ being minimal.
This shows that the assumption $f_s<p^k$ was wrong, so we have $f_s=p^k$.
Note that $|s^{(p^k)}|=0$ implies $p|\lambda(s,p^k)$ and $p|\Tilde{\lambda}(s,p^k)$ since if $|s^{(p^k)}|=0$, then $s^{(p^k)}$ is a union of vanishing sums of the form $(z^a,z^{a+p^{k-1}},\dots,z^{a+p^{k-1}(p-1)})$ for some $a\in\mathbb{Z}$ with $p\nmid a$ which shows the last claim.
\item First assume that $|s^{(p^k)}|=0$.
In this case $s$ is a sum of $p^{k-1}$-th roots of unity, so $f_s<p^k$ and therefore
\[[\mathbb{Q}_{p^k}:\mathbb{Q}(s)]=[\mathbb{Q}_{p^k}:\mathbb{Q}_{p^{k-1}}]\cdot[\mathbb{Q}_{p^{k-1}}:\mathbb{Q}(s)]=p\cdot[\mathbb{Q}_{p^{k-1}}:\mathbb{Q}(s)],\]
so in particular $p|[\mathbb{Q}_{p^k}:\mathbb{Q}(s)]$.

For the other direction we can assume by part (a) that $s$ is minimal since $f_s$ and $\mathbb{Q}(s)$ only depend on $|s|$.
Let $\sigma\in\Gal(\mathbb{Q}_{p^k}/\mathbb{Q}(s))$ be a Galois automorphism.
Now $\sigma$ fixes $|s|$ and since it is a Galois automorphism, it maps roots of unity to roots of unity which leads to the following equation:
\[z^{a_1}+\dots+z^{a_l}=\sigma(z^{a_1})+\dots+\sigma(z^{a_l}).\]
So $s=(z^{a_1},\dots,z^{a_l})$ and $s':=(\sigma(z^{a_1}),\dots,\sigma(z^{a_l}))$ are sums of $p^k$-th roots of unity with the same value and the same length.
Since $s$ is minimal, it follows from (a) that $s'$ is minimal as well and that $s=s'$ which shows that $\sigma$ permutes the entries of $s$.
This implies that $s$ is a union of $\sigma$-orbits and that $s^{(p^k)}$ is a union of $\sigma$-orbits which consist of primitive $p^k$-th roots of unity.
Now suppose that $o(\sigma)=p$.
Then there is some $b\in\{1,\dots,p-1\}$ such that for all $a\in\mathbb{Z}$ we have
\[(*)\ \sigma(z^a)=z^{a\cdot(b\cdot p^{k-1}+1)}.\]
Indeed, we have
\[(b\cdot p^{k-1}+1)^p=\sum\limits_{i=0}^p\binom{p}{i}\cdot(b\cdot p^{k-1})^i\cdot1^{p-i}=\sum\limits_{i=0}^p\binom{p}{i}\cdot(b\cdot p^{k-1})^i\]
and only the summand for $i=0$ is not divisible by $p^k$, so
\[(b\cdot p^{k-1}+1)^p\equiv\binom{p}{0}\cdot(b\cdot p^{k-1})^0\equiv1\mod p^k,\]
which shows that if $\sigma$ satisfies $(*)$, then it is indeed an automorphism of order $p$.
On the other hand, since $\varphi(p^k)=(p-1)\cdot p^{k-1}$, there are exactly $p-1$ choices for $\sigma$ which correspond to the choice of $b$.
Further, if $c\in\{1,\dots,p-1\}$, then
\[(b\cdot p^{k-1}+1)^c=\sum\limits_{i=0}^c\binom{c}{i}\cdot(b\cdot p^{k-1})^i\equiv c\cdot b\cdot p^{k-1}+1\mod p^k,\]
which shows that independently of $b$, each $\sigma$-orbit which consists of primitive $p^k$-th roots of unity is of the form $(z^a,z^{a+p^{k-1}},\dots,z^{a+p^{k-1}(p-1)})$ for some $a\in\mathbb{Z}$ with $p\nmid a$, so it is a vanishing sum.
It follows that $|s^{(p^k)}|=0$ which proves the first claim.
For the second claim, note that, like in (b), $|s^{(p^k)}|=0$ implies $p|\lambda(s,p^k)$ and $p|\Tilde{\lambda}(s,p^k)$.\qedhere
\end{compactenum}
\end{proof}
\end{theorem}

With the help of Theorem 2.5 we can now prove that Property 1.3 is equivalent to a seemingly weaker statement as follows.
\begin{corollary}
Let $G$ be a finite group, let $p$ be an odd prime, let $\chi\in\Irr_{p'}(G)$ and let $f$ be the conductor of $\chi$.
Assume that $f_p\geq p^2$, which can be done by Lemma 2.1.
Then Property 1.3 holds for $G,\ \chi$ and $p$ if and only if there is some $p$-element $g\in G$ with $(f_{\chi})_p=(f_{\chi(g)})_p$.
\begin{proof}
It is clear that if Property 1.3 holds for $\chi$ with $g\in G$ being the corresponding $p$-element, then $f_p=(f_{\chi(g)})_p$, since
\[[\mathbb{Q}_{f_p}:\mathbb{Q}(\chi(g))]=[\mathbb{Q}_{f_p}:\mathbb{Q}_{(f_{\chi(g)})_p}]\cdot[\mathbb{Q}_{(f_{\chi(g)})_p}:\mathbb{Q}(\chi(g))]\]
and since $p|[\mathbb{Q}_{f_p}:\mathbb{Q}_{(f_{\chi(g)})_p}]$ if $f_p\ne (f_{\chi(g)})_p$.

For the converse, let $g\in G$ be a $p$-element with $f_p=(f_{\chi(g)})_p$ and write $p^a$ for $f_p$.
Further, let $s$ be a sum of $p^k$-th roots of unity for some $k\in\mathbb{Z}_{\geq2}$ with $|s|=\chi(g)$.
Then by Theorem 2.5 (b), we have $|s^{(p^b)}|=0$ for all $b\in\mathbb{Z}_{\geq1}$ with $b>a$ and $|s^{(p^a)}|\ne0$.
It follows that $p\nmid[\mathbb{Q}_{p^a}:\mathbb{Q}(s)]$ by Theorem 2.5 (c), which shows that Property 1.3 holds for $\chi$.
\end{proof}
\end{corollary}

Let again $k\in\mathbb{Z}_{\geq2}$ and let $c\in\mathbb{C}$ such that there is a sum of $p^k$-th roots of unity $s$ with $|s|=c$.
For $m\in\mathbb{Z}_{\geq2}$ we set
\[\lambda(c,p^m):=\lambda(s,p^m),\ \Tilde{\lambda}(c,p^m):=\Tilde{\lambda}(s,p^m)\]
considered as elements of $\mathbb{Z}/p\mathbb{Z}$.
This is well-defined by Theorem 2.5 (a)+(b).

\section{Series of finite groups with known character table}

We start by investigating some infinite series of finite quasi-simple groups for which the character table is known.
By using the theory we developed, we get that these groups have Property 1.3 as shown by the following theorem.

\begin{theorem}
Let $G$ be a finite group, let $q$ be a prime power and let $p$ be an odd prime.
Further, let $\chi\in\Irr(G)$ and set $S:=G/\Z(G)$.
Then Property 1.3 holds for $G,\ \chi$ and $p$ in all of the following cases.
\begin{itemize}
\item $S\in\{\PSL_2(q),\PSL_3(q),\PSU_3(q),\PSp_4(q),G_2(q),{^3}D_4(q)\}$ and $q$ is any prime power
\item $S\in\{\PSp_6(2^n),{^2}B_2(2^{2n+1}),{^2}G_2(3^{2n+1})\}$ for some $n\in\mathbb{Z}_{\geq1}$
\item $S={^2}F_4(2)'$
\item $S={^2}F_4(q)$, $q=2^{2n+1}$ for some $n\in\mathbb{Z}_{\geq1}$ and $\chi$ does not lie in one of the ten character families not given in \cite{CHEVIE}
\end{itemize}
\begin{proof}
First note that if $G$ has Property 1.3, then any factor group $H$ of $G$ also has Property 1.3 since for every character $\chi\in\Irr(H)$ its inflation $\Inf_H^G(\chi)$ is an irreducible character of $G$ with the same degree and the same field of values as $\chi$.
From this we see that it is enough to consider the case where $G$ is a central extension of $S$ by the Schur multiplier of $S$.

A list of exceptional, i.e., non generic Schur multipliers of the groups we want to consider is given in \cite[Table 24.3]{MaTe}.
If $G$ is the Tits group ${^2}F_4(2)'$ or a central extension of $S$ by an exceptional Schur multiplier, then Property 1.3 can be checked directly using \cite{GAP}, so we can assume that $G$ is a central extension of $S$ by a generic Schur multiplier.

The proof works as follows.
Let $\chi\in\Irr(G)$ and let $f$ be its conductor.
We can assume that $p^2|f$ by Lemma 2.1 and that $p\nmid\chi(1)$.
In all cases we want to consider, there is a unique irreducible polynomial $Q$ in $q$ with $p^2|Q$ and such that for all $g\in G$ the polynomial $Q$ divides the order of every root of unity of order divisible by $p$ which appears in $\chi(g)$ as polynomials in $q$.
Note that the character tables we want to consider are generic character tables of finite groups of Lie type, so they depend on the prime power $q$.
We try to find a $p$-power $f'_p\in\mathbb{Z}_{\geq1}$ such that for the $p$-part of $f_{\chi(h)}$ we have $(f_{\chi(h)})_p\leq f'_p$ for all $h\in G$, so $f_p\leq f'_p$ and such that there is a $p$-element $g\in G$ with $\lambda(\chi(g),f'_p)\ne0$ which is then of order $Q_p$

If such a $p$-element $g$ exists, then $f_p=f'_p$ by Theorem 2.5 (b) and then Property 1.3 holds by Theorem 2.5 (c).
Note that if $s$ is a sum of $f'_p$-th roots of unity with $|s|=\chi(g)$, then $\lambda(\chi(g),f'_p)\ne0$ and $\Tilde{\lambda}(\chi(g),f'_p)\ne0$ each imply $|s^{(f'_p)}|\ne0$ as in the proof of Theorem 2.5 (b).
The tables in Section 5 give such a $p$-power $f'_p$ and a suitable $p$-element $g$ for all cases where $p^2|f$ and $p\nmid\chi(1)$.
For each table we use the notation from the reference in which the character table of the respective family of groups is given.

Note that in \cite{CHEVIE} there is a unique pair of numbers $[a,b]$ associated to each irreducible character and each conjugacy class of $G$.
We write $\chi_{a,b}$ for the character corresponding to $[a,b]$ and $C_a$ for the conjugacy class corresponding to $[a,0]$ whenever \cite{CHEVIE} is given as the reference for the character table.
The conjugacy classes with $b=0$ are the semisimple conjugacy classes and clearly any $p$-element for $p\nmid q$ has to lie in a semisimple conjugacy class.

If $p$ does not divide $\Tilde{\lambda}(\chi(g),f'_p)$, then Property 1.3 holds by Theorem 2.5 (b)+(c), so it is left to consider the cases in which $p$ divides $\Tilde{\lambda}(\chi(g),f'_p)$ which is done in the following.

If $p=3$, then $p|(q-1)$ if and only if $p|(q^2+q+1)$ and $p|(q+1)$ if and only if $p|(q^2-q+1)$ which implies that $p|\chi(1)$ in all cases of Table 3, 4, 7, 8, 9, 10 and 11 in which $3$ might divide $\Tilde{\lambda}(\chi(g),f'_p)$ so Property 1.3 does not make a statement in these cases.

If $G={^2}G_2(q)$ and $p=3$, then $p\nmid Q$ in all cases of Table 14 since $q=3^{2n+1}$ and if $p=5$, then $p\nmid(q+1)$, since
\[q=3^{2n+1}\equiv i\mod5\ \text{ for some }i\in\{2,3\},\]
so these cases are not possible.

If $G={^2}F_4(q)$ and $p=3$, then $p$ does not divide any of $q-1,\ Q_1,\ Q_2$ since
\[q=2^{2n+1}\equiv2\mod3\]
and therefore
\[q-1\equiv Q_1\equiv Q_2\equiv1\mod3.\]
Note that $3|(q^2-q+1)$, so $p$ divides the degree of the characters $\chi_{15,0}(k,l)$ of Table 13 and then Property 1.3 does not make a statement.

If $G={^2}F_4(q)$ and $p=5$, then $p$ does not divide any of $q+1,\ Q_1,\ Q_2$ since
\[q=2^{2n+1}\equiv i\mod5\ \text{ for some }i\in\{2,3\}\]
and therefore
\begin{align*}
q+1&\equiv i\mod5\ \text{ for some }i\in\{3,4\}\text{ and}\\
Q_j&\equiv i\mod5\ \text{ for some }i\in\{1,2,3,4\}\text{ and all }j\in\{1,2\}.
\end{align*}
If $G={^2}F_4(q)$ and $p=7$, then also $p\nmid(q+1)$ since
\[q=2^{2n+1}\equiv i\mod7\ \text{ for some }i\in\{1,2,4\}\]
and therefore
\[q+1\equiv i\mod7\ \text{ for some }i\in\{2,3,5\}\]
so these cases are not possible.

It remains to consider the case $G={^2}F_4(q)$ and $p=11$ for the character family $\chi=\chi_{15,0}(k,l)$ of Table 13.
Here $\chi(g)$ is, up to multiplication by a rational number, a sum of $24$ roots of unity.
We set $r:=2^n$.
Now all exponents of these roots of unity are multiples of the parameter $a$ of the conjugacy class which is not divisible by $p$ and since we are only interested in the number of these exponents which are divisible by $p$ we can leave out the factor $a$ in the exponents.
Then the exponents of $12$ of these roots of unity are
\[\{k+l,k-l,2rk,2rl,k+rk+rl,-k-rk+rl,k-rk+rl,\]
\[-k+rk+rl,rk-l-rl,rk+l+rl,rk+l-rl,+rk-l+rl\}\]
and the other $12$ are their respective negatives.
It is enough to show that it is not possible to choose $22$ of the $24$
roots of unity such that they build two minimal vanishing sums of length $11$.
Now suppose that it is possible and let $E$ be the set consisting of the exponents of these $22$ roots of unity.
Then the two exponents not contained in $E$ are divisible by $p$ and they add up to $0$.
We can write
\[E=E_1\cup E_2\]
for the two subsets $E_1$ and $E_2$ of $E$ with $|E_1|=|E_2|=11$ whose elements are the exponents of roots of unity which build a minimal vanishing sum.
It follows that if $a,b\in E_1$, then $p|(a-b)$ and the same holds for $E_2$.
Further, if $a\in E_1$, then $p\nmid a$ and $p\nmid-a$, so $-a\in E$ and since $a$ and $-a$ cannot correspond to the same minimal vanishing sum, we have $-a\in E_2$.
Analogously, if $a\in E_2$, then $-a\in E_1$.

Now suppose that $2rk,-2rk\notin E_1$.
Then we have $p|k$ since $p\nmid(2r)$ and we have that $E_1$ contains one of $\pm(k+rk+rl)$ and one of $\pm2rl$.
Let $a,\ b$ be these two elements of $E_1$.
Then $p|(a-b)$ and
\[a-b\in\{\pm(k+rk-rl),\pm(k+rk+3rl)\}.\]
In either case we get $p|l$ since $p\nmid r,\ p\nmid3r$, but then $p$ divides every element of $E$ which is a contradiction to $\lambda(\chi(g),f'_p)\geq1$.
Analogously, it can be shown that one of $\pm2rl$ lies in $E_1$.

Now assume that $2rk,2rl\in E_1$.
Then $p|(2rk-2rl)$, so $p|(k-l)$, which shows that $\pm(k-l)$ are the two exponents not contained in $E$.
Further, $E_1$ contains one of $\pm(k+rk+rl)$ and one of $\pm(-k+rk+rl)$, so $p$ divides the difference of these two which lies in
\[\{\pm2k,\pm(2rk+2rl)\}.\]
It follows that either $p|k$ or $p|(k+l)$ which contradicts $2rk\in E$ or $k+l\in E$ respectively.

If instead $2rk,-2rl\in E_1$, then $p|(2rk+2rl)$, which implies $p|(k+l)$ and then $\pm(k+l)$ are the two exponents not contained in $E$.
Here $E_1$ also contains one of $\pm(k+rk-rl)$ and one of $\pm(-k+rk-rl)$, so $p$ divides the difference of these two which lies in
\[\{\pm2k,\pm(2rk-2rl)\}.\]
Then we have that either $p|k$ or $p|(k-l)$ which contradicts $2rk\in E$ or $k-l\in E$ respectively.

The two cases where $-2rk,2rl\in E_1$ and $-2rk,-2rl\in E_1$ are analogous to the other two cases.
This finishes the proof.
\end{proof}
\end{theorem}

\section{The finite general and special linear groups}

In this section we continue our investigation of the finite linear groups.
We develop some results on Property 1.3 for the finite general and special linear groups in arbitrary dimension.

The irreducible complex characters of the finite general linear groups have been parameterised in \cite{Gre} and we start by summarising the results developed there.
See also \cite{Gup}, where the work of Green has been used to compute the character table of $\GL_5(q)$ for all prime powers $q$.

A \textbf{partition} $\rho=(a_1,\dots,a_l)$ of some $n\in\mathbb{Z}_{\geq0}$ is a finite non-increasing sequence of numbers $a_i\in\mathbb{Z}_{\geq1}$ such that $|\rho|:=a_1+\dots+a_l=n$.
Its entries $a_1,\dots,a_l$ are called \textbf{parts} of $\rho$.
Note that the only partition of $0$ is the empty partition $()$.
Now if $\mathfrak{f}(t)=t^d+c_{d-1}t^{d-1}+\dots+c_0$ is a polynomial with coefficients in the finite field $\mathbb{F}_q$, then the \textbf{companion matrix} of $\mathfrak{f}$ is defined as
\[U(\mathfrak{f}):=\begin{pmatrix}&1&&&\\&&\ddots&&\\&&&1&\\-c_0&-c_1&-c_2&\dots&-c_{d-1}\end{pmatrix}\in\M^{d\times d}(\mathbb{F}_q).\]
Further, for $m\in\mathbb{Z}_{\geq1}$ we set
\[U_m(\mathfrak{f}):=\begin{pmatrix}U(\mathfrak{f})&I_d&&\\&\ddots&\ddots&\\&&\ddots&I_d\\&&&U(\mathfrak{f})\end{pmatrix}\in\M^{md\times md}(\mathbb{F}_q)\]
with $m$ diagonal blocks $U(\mathfrak{f})$, where $I_n$ is the identity matrix in dimension $n$ and if $\rho=(a_1,\dots,a_l)$ is a partition of some $b\in\mathbb{Z}_{\geq0}$, then
\[U_\rho(\mathfrak{f}):=\diag(U_{a_1}(\mathfrak{f}),\dots,U_{a_l}(\mathfrak{f}))\in\M^{bd\times bd}(\mathbb{F}_q).\]
Now if $A\in\GL_n(q)$ for some $n\in\mathbb{Z}_{\geq1}$ with characteristic polynomial $\mathfrak{f}_1^{n_1}\dots\mathfrak{f}_k^{n_k}$, where $\mathfrak{f}_1,\dots,\mathfrak{f}_k$ are distinct irreducible polynomials over $\mathbb{F}_q$, then, as in \cite{Gre}, $A$ is conjugate to a matrix $\diag(U_{\rho_1}(\mathfrak{f}_1),\dots,U_{\rho_k}(\mathfrak{f}_k))$, where $\rho_i$ is a partition of $n_i$ for all $i\in\{1,\dots,k\}$ and
\[\sum\limits_{i=1}^k|\rho_i|\cdot\deg(\mathfrak{f}_i)=n.\]
Conversely, let $\mathfrak{F}$ be the set of all irreducible polynomials over $\mathbb{F}_q$ and let $P$ be the set of all partitions.
Then by \cite[Lemma 1.1]{Gre} a conjugacy class $c$ of $\GL_n(q)$ can be specified by the function $\nu_c:\mathfrak{F}\rightarrow P$ which maps an irreducible polynomial $\mathfrak{f}$ to the partition with which $\mathfrak{f}$ is associated in the canonical form of a matrix of $c$ and then
\[\sum\limits_{\mathfrak{f}\in F}|\nu_c(\mathfrak{f})|\cdot\deg(\mathfrak{f})=n.\]
Thus we get the following Proposition.

\begin{proposition}
Let $\mathfrak{F}_n$ be the set of all maps $\nu:\mathfrak{F}\rightarrow P$ such that
\[\sum\limits_{\mathfrak{f}\in\mathfrak{F}}|\nu(\mathfrak{f})|\cdot\deg(\mathfrak{f})=n.\]
Then the map $\Phi:\Cl(\GL_n(q))\rightarrow\mathfrak{F}_n,\ c\mapsto\nu_c$ is a bijection.
\end{proposition}

We denote a conjugacy class $c\in\Cl(\GL_n(q))$ by the symbol $(\dots \mathfrak{f}^{\nu_c(\mathfrak{f})}\dots)$.

There is a similar parameterisation for the set of irreducible complex characters of $\GL_n(q)$ for which we need the notion of a \textbf{$d$-simplex}.
Let $\mathfrak{f}$ be an irreducible polynomial of degree $d$ over $\mathbb{F}_q$ and let $\omega$ be a generator of $\mathbb{F}_{q^d}^\times$.
Then a root of $\mathfrak{f}$ is of the form $\omega^k$ for some $k\in\mathbb{Z}_{\geq0}$ and the set $\{\omega^k,\omega^{kq},\dots,\omega^{kq^{d-1}}\}$ consists of all the roots of $\mathfrak{f}$ by \cite[Section 7]{Gre}.
In this case the set $\{k,kq,\dots,kq^{d-1}\}$ is called $d$-simplex with $k,kq,\dots,kq^{d-1}$ as its roots and $d$ as its degree.
We denote the set of all simpleces by $\mathfrak{G}$.

Now by \cite[Theorem 14]{Gre} every $\chi\in\Irr(\GL_n(q))$ is characterised by a map $\nu_\chi:\mathfrak{G}\rightarrow P$ with
\[\sum\limits_{\mathfrak{g}\in\mathfrak{G}}|\nu_\chi(\mathfrak{g})|\cdot\deg(\mathfrak{g})=n\]
and the following Proposition, which is one of the claims of \cite[Theorem 14]{Gre}, holds.

\begin{proposition}
Let $\mathfrak{G}_n$ be the set of all maps $\nu:\mathfrak{G}\rightarrow P$ such that
\[\sum\limits_{\mathfrak{g}\in\mathfrak{G}}|\nu(\mathfrak{g})|\cdot\deg(\mathfrak{g})=n.\]
Then the map $\Psi:\Irr(\GL_n(q))\rightarrow\mathfrak{G}_n,\ \chi\mapsto\nu_\chi$ is a bijection.
\end{proposition}

We denote an irreducible complex character $\chi\in\Irr(\GL_n(q))$ by the symbol $(\dots\mathfrak{g}^{\nu_\chi(\mathfrak{g})}\dots)$.

Now let $\rho$ be a partition.
A \textbf{set of $\rho$-variables} is a formal set
\[X^\rho=\{x_{1,1},\dots,x_{1,r_1},x_{2,1},\dots,x_{2,r_2},\dots\},\]
where $r_i$ is the number of parts of $\rho$ which are equal to $i$ and the elements of $X^\rho$ are called \textbf{$\rho$-variables}.
We define the degree of a $\rho$-variable $x_{d,i}$ as $\deg(x_{d,i}):=d$.
Now to each $\rho$-variable $x_{d,i}$ of degree $d$ there correspond $d$ variables called the \textbf{roots} of $x_{d,i}$ or \textbf{$\rho$-roots}, which are denoted as $\xi_{d,i},\xi_{d,i}^q,\dots,\xi_{d,i}^{q^{d-1}}$.

A \textbf{substitution of $X^\rho$ into $\mathfrak{F}$} is a map $\alpha:X^\rho\rightarrow\mathfrak{F}$ such that
\[\deg(\alpha(x))|\deg(x)\]
for all $x\in X^\rho$ and for each $\rho$-root
$\xi_{d,i}$ we choose any root $\gamma_{d,i}$ of $\alpha(x_{d,i})$ and set $\alpha(\xi_{d,i}^{q^j}):=\gamma_{d,i}^{q^j}$.
Two substitutions $\alpha,\ \beta$ are called \textbf{equivalent} if there is a degree preserving permutation $\sigma$ of $X^\rho$ such that $\alpha=\beta\circ\sigma$.
Further, a class of equivalent substitutions is called \textbf{mode of substitution}.

Now let $\alpha$ be a substitution of the set of $\rho$-variables $X^\rho$ into $\mathfrak{F}$ and let $\mathfrak{f}\in\mathfrak{F}$ be of degree $d$.
For each positive integer $i$ we set
\[r_i(\alpha,\mathfrak{f}):=|\{x\in X^\rho\ |\ \deg(x)=i\cdot d\text{ and }\alpha(x)=\mathfrak{f}\}|.\]
Then we define the partition $\rho(\alpha,\mathfrak{f}):=\{1^{r_1(\alpha,\mathfrak{f})}2^{r_2(\alpha,\mathfrak{f})}\dots\}$ which means that $\rho(\alpha,\mathfrak{f})$ has $r_i(\alpha,\mathfrak{f})$ parts which are equal to $i$ for every $i\in\mathbb{Z}_{\geq1}$.
It is easily seen that $\rho(\alpha,\mathfrak{f})$ only depends on the mode of substitution $m$ of $\alpha$ so we can define $\rho(m,\mathfrak{f}):=\rho(\alpha,\mathfrak{f})$ without ambiguity.

Finally, if $c=(\dots\mathfrak{f}^{\nu_c(\mathfrak{f})}\dots)$ is a conjugacy class of $\GL_n(q)$, then a substitution $\alpha$ of $X^\rho$ into $\mathfrak{F}$ such that $|\rho(\alpha,\mathfrak{f})|=|\nu_c(\mathfrak{f})|$ for all $\mathfrak{f}\in\mathfrak{F}$ is called \textbf{substitution of $X^\rho$ into the class $c$}.

The notions of a \textbf{substitution of $X^\rho$ into $\mathfrak{G}$} and a \textbf{substitution of $X^\rho$ into the character $\chi$} are defined analogously.

Before we are able to prove the main theorem for the finite general linear groups in arbitrary dimension we need to define \textbf{cyclotomic polynomials} and we give a lemma on the degrees of the irreducible complex characters of $\GL_n(q)$.

In the complex numbers all primitive $n$-th roots of unity have a common minimal polynomial in $\mathbb{Q}[x]$, the $n$-th cyclotomic polynomial $\Phi_n(x)$ (see \cite[Section 5]{Con}) and we have
\[x^n-1=\prod\limits_{d|n}\Phi_d(x),\]
so $\Phi_i(q)|(q^j-1)$ as polynomials in $q$ if and only if $i|j$.
Note that if $\Phi_i(q)|(q^j-1)$ as polynomials in $q$, then also $\Phi_i(k)|(k^j-1)$ for all $k\in\mathbb{Z}$.

For the following lemma, we need another notation.
If $i\in\mathbb{Z}$ and $j\in\mathbb{Z}_{\geq1}$, we set
\[\lfloor i\rfloor_j:=\max\{k\in\mathbb{Z}\ |\ k\leq i,\ j|k\},\]
which is $i$ rounded down to the next integer divisible by $j$.

\begin{lemma}
Let $\chi=(\dots\mathfrak{g}^{\nu_\chi(\mathfrak{g})}\dots)\in\Irr(\GL_n(q))$ for some $n\in\mathbb{Z}_{\geq1}$ and some prime power $q\ne1$ and let $m\in\mathbb{Z}_{\geq1}$.
If
\[\sum\limits_{\mathfrak{g}\in\mathfrak{G}}\bigg(|\nu_\chi(\mathfrak{g})|\cdot\deg(\mathfrak{g})-\Big\lfloor|\nu_\chi(\mathfrak{g})|\cdot\gcd(\deg(\mathfrak{g}),m)\Big\rfloor_m\bigg)\geq a\cdot m,\]
then $\Phi_m(q)^a|\chi(1)$ as polynomials in $q$ for all $a\in\mathbb{Z}_{\geq0}$.
\begin{proof}
By \cite[Theorem 14]{Gre} the degree $\chi(1)$ of $\chi$ is given by
\[\Psi_n(q)\prod\limits_{\mathfrak{g}\in\mathfrak{G}}(-1)^{|\nu_\chi(\mathfrak{g})|}\{\nu_\chi(\mathfrak{g}):q^{\deg(\mathfrak{g})}\},\]
where
\[\Psi_i(q):=\prod\limits_{j=1}^i(q^j-1)\text{ for }i\in\mathbb{Z}_{\geq0}\]
and for any partition $\rho=(a_1,\dots,a_l)$ we have
\[\{\rho:q\}:=q^{a_2+2a_3+\dots}\prod\limits_{1\leq r<s\leq l}(1-q^{a_r-a_s-r+s})\bigg/\prod\limits_{r=1}^l(-1)^{a_r+l-r}\Psi_{a_r+l-r}(q).\]
Further, if $Q\in\mathbb{Q}(q)$, the field of rational functions in $q$ over $\mathbb{Q}$, we set
\[F_m(Q):=\max\bigg\{i\in\mathbb{Z}\ |\ \frac{Q}{\Phi_m(q)^i}\in\mathbb{Z}[q]\bigg\}.\]
Since $\Phi_i(q)|(q^j-1)$ as polynomials in $q$ if and only if $i|j$, we get that
\[F_m(\{\rho:q^k\})=F_m(\{\rho:q^{\gcd(k,m)}\}).\]
Now we choose for each simplex $\mathfrak{g}\in\mathfrak{G}$ with $|\nu_\chi(\mathfrak{g})|\geq1$ a simplex $\mathfrak{g}'\in\mathfrak{G}$ with
\[\deg(\mathfrak{g}')=\gcd(\deg(\mathfrak{g}),m)\]
and such that all the $\mathfrak{g}'$ are different.
For a small prime power $q$ there might not be enough different $\mathfrak{g}'$, but
it is still possible to choose the $\mathfrak{g}'$ in this way since $q$ is considered to be a formal variable here for which the prime power is inserted later.
Let $\chi'$ be the irreducible complex character of the possibly different group $\GL_{n'}(q)$ defined by
\[\nu_{\chi'}(\mathfrak{g}'):=\nu_{\chi}(\mathfrak{g})\]
for all $\mathfrak{g}\in\mathfrak{G}$ with $|\nu_\chi(\mathfrak{g})|\geq1$ and $()$, the empty partition, otherwise.
Then
\[\prod\limits_{\mathfrak{g}\in\mathfrak{G}}(-1)^{|\nu_{\chi'}(\mathfrak{g})|}\{\nu_{\chi'}(\mathfrak{g}):q^{\deg(\mathfrak{g})}\}=\prod\limits_{\mathfrak{g}\in\mathfrak{G}}(-1)^{|\nu_\chi(\mathfrak{g})|}\{\nu_\chi(\mathfrak{g}):q^{\deg(\mathfrak{g})}\}\]
and we have
\[n'=\sum\limits_{\mathfrak{g}\in\mathfrak{G}}|\nu_\chi(\mathfrak{g})|\cdot\gcd(\deg(\mathfrak{g}),m).\]
Now let $\mathfrak{g}\in\mathfrak{G}$ be a simplex of degree $d$ with $d|m$ and let $\chi_\mathfrak{g}\in\GL_{d\cdot|\nu_\chi(\mathfrak{g})|}(q)$ and $\overline{\chi_\mathfrak{g}}\in\Irr(\GL_{d\cdot\lfloor|\nu_\chi(\mathfrak{g})|\rfloor_{m/d}}(q))$ be defined by
\[\nu_{\chi_\mathfrak{g}}(\mathfrak{g})=\nu_\chi(\mathfrak{g}),\ \nu_{\overline{\chi_\mathfrak{g}}}(\mathfrak{g})=\{\big\lfloor|\nu_\chi(\mathfrak{g})|\big\rfloor_{m/d}^1\}\]
and $()$ otherwise.
Further, by \cite[Remark 2.3.27]{GeMa} the degree of any character of $\GL_n(q)$ is a polynomial in $q$ so $F_m(\chi_\mathfrak{g}(1))\geq0$.
By using this we get
\[0\leq F_m(\chi_\mathfrak{g}(1))=\left\lfloor\frac{d\cdot|\nu_\chi(\mathfrak{g})|}{m}\right\rfloor+F_m(\{\nu_\chi(\mathfrak{g}):q^d\})\]
and therefore
\[F_m(\{\nu_\chi(\mathfrak{g}):q^d\})\geq-\left\lfloor\frac{d\cdot|\nu_\chi(\mathfrak{g})|}{m}\right\rfloor.\]
On the other hand we have
\[F_m(\{\nu_{\overline{\chi_\mathfrak{g}}}(\mathfrak{g}):q^d\})=-\left\lfloor\frac{d\cdot|\nu_{\overline{\chi_\mathfrak{g}}}(\mathfrak{g})|}{m}\right\rfloor=-\left\lfloor\frac{d\cdot\big\lfloor|\nu_\chi(\mathfrak{g})|\big\rfloor_{m/d}}{m}\right\rfloor=-\left\lfloor\frac{d\cdot|\nu_\chi(\mathfrak{g})|}{m}\right\rfloor\]
so for all $\mathfrak{g}\in\mathfrak{G}$ we have
\[F_m(\{\nu_\chi(\mathfrak{g}):q^{\deg(\mathfrak{g})}\})\geq F_m(\{\nu_{\overline{\chi_\mathfrak{g}}}(\mathfrak{g}):q^{\deg(\mathfrak{g})}\}).\]
We define $\chi'_\mathfrak{g}$ and $\overline{\chi'_\mathfrak{g}}$ analogously.
Now let $\overline{\chi}\in\Irr(\GL_{\overline{n}}(q))$ be defined by
\[\nu_{\overline{\chi}}(\mathfrak{g})=\nu_{\overline{\chi'_\mathfrak{g}}}(\mathfrak{g})\]
for all $\mathfrak{g}\in\mathfrak{G}$, so
\[\overline{n}=\sum\limits_{\mathfrak{g}\in\mathfrak{G}}\bigg(\Big\lfloor|\nu_\chi(\mathfrak{g})|\cdot\gcd(\deg(\mathfrak{g}),m)\Big\rfloor_m\bigg).\]
Finally, we can do the following calculation
\begin{align*}
F_m(\chi(1))&=F_m(\Psi_n(q))+\sum\limits_{\mathfrak{g}\in\mathfrak{G}}F_m((-1)^{|\nu_\chi(\mathfrak{g})|}\{\nu_\chi(\mathfrak{g}):q^{\deg(\mathfrak{g})}\})\\
&=F_m(\Psi_n(q))+\sum\limits_{\mathfrak{g}\in\mathfrak{G}}F_m((-1)^{|\nu_{\chi'}(\mathfrak{g})|}\{\nu_{\chi'}(\mathfrak{g}):q^{\deg(\mathfrak{g})}\})\\
&\geq F_m(\Psi_n(q))+\sum\limits_{\mathfrak{g}\in\mathfrak{G}}F_m((-1)^{|\nu_{\overline{\chi'_\mathfrak{g}}}(\mathfrak{g})|}\{\nu_{\overline{\chi'_\mathfrak{g}}}(\mathfrak{g}):q^{\deg(\mathfrak{g})}\})\\
&=F_m(\Psi_n(q))-F_m(\Psi_{\overline{n}}(q))+F_m(\overline{\chi}(1))\\
&\geq F_m(\Psi_n(q))-F_m(\Psi_{\overline{n}}(q))\\
&=F_m(\Psi_n(q)/\Psi_{\overline{n}}(q))\\
&=F_m((q^{\overline{n}+1}-1)\cdot(q^{\overline{n}+2}-1)\dots(q^n-1)).
\end{align*}
This is at least $a$ if $n-\overline{n}\geq a\cdot m$ since $\Phi_i(q)|(q^j-1)$ as polynomials in $q$ if and only if $i|j$.
The condition $n-\overline{n}\geq a\cdot m$ is equivalent to the assumption
\[\sum\limits_{\mathfrak{g}\in\mathfrak{G}}\bigg(|\nu_\chi(\mathfrak{g})|\cdot\deg(\mathfrak{g})-\Big\lfloor|\nu_\chi(\mathfrak{g})|\cdot\gcd(\deg(\mathfrak{g}),m)\Big\rfloor_m\bigg)\geq a\cdot m,\]
which finishes the proof.
\end{proof}
\end{lemma}

With this lemma we can now prove the main theorem about the finite general linear groups in arbitrary dimension.
In the proof of the following Theorem we choose the $p$-element $g$ needed for Property 1.3 such that $\det(g)=1$ whenever possible.
This helps us later to prove that Property 1.3 also holds for the finite special linear groups under certain assumptions.

\begin{theorem}
Let $n\in\mathbb{Z}_{\geq1}$ and let $q$ be a prime power.
Further, let $p$ be a prime with $p|(q-1)$ and $p>n(n-1)$.
Then Property 1.3 holds for $\GL_n(q)$ and $p$.
\begin{proof}
Let $\chi\in\Irr(\GL_n(q))$, let $f$ be the conductor of $\chi$ and let $c\in\Cl(\GL_n(q))$.
If $m_\chi$ is a mode of substitution of $X^\rho$ into $\chi$ for some partition $\rho$, we write $\mu(m_\chi,\chi)$ for the value $\chi(m_\chi,\chi)$ defined in \cite{Gre} to avoid confusion with the character $\chi$.
Then by \cite[Theorem 14]{Gre} the value of $\chi$ on the class $c$ is given by
\begin{align*}
\chi(c)=\ &(-1)^{n-\sum\limits_{\mathfrak{g}\in\mathfrak{G}}|\nu_\chi(\mathfrak{g})|}\sum\limits_\rho\sum\limits_{m_\chi}\mu(m_\chi,\chi)\sum\limits_{m_c}Q(m_c,c)\cdot\\
&\prod\limits_d\sum\limits_{\sigma\in\S_{r_d}}S_d(h_{d1}m_\chi:\xi_{d\sigma(1)})\cdots S_d(h_{dr_d}m_\chi:\xi_{d\sigma(r_d)}).
\end{align*}
Here $Q(m_c,c)$ is defined as in \cite{Gre},
\[S_d(k:\xi)=\theta^k(\xi)+\theta^{qk}(\xi)+\dots+\theta^{q^{d-1}k}(\xi),\ k\in\mathbb{Z}_{\geq1},\]
where $\theta$ is a generator of the character group of $\mathbb{F}_{q^{n!}}^\times$ and $\xi$ is a root of an irreducible polynomial over $\mathbb{F}_q$.
Further,
\[h_{di}m=c_{di}\frac{q^d-1}{q^{s_{di}}-1},\]
where if $\alpha$ is a substitution which belongs to the mode of substitution $m$, then $c_{di}$ is a root of the simplex $\alpha(x_{di})$ and $s_{di}$ is the degree of $\alpha(x_{di})$.
The sums and the product are taken over all partitions $\rho=\{1^{r_1}2^{r_2}\dots\}$ of $n$, all modes of substitutions $m_\chi$ of $X^\rho$ into $\chi$, all modes of substitutions $m_c$ of $X^\rho$ into $c$, all $d\in\mathbb{Z}_{\geq1}$ with $r_d\geq1$ and all permutations $\sigma\in\S_{r_d}$, where $\S_n$ is the symmetric group on the set $\{1,\dots,n\}$, respectively.

Now suppose that there is a simplex $\mathfrak{g}_0\in\mathfrak{G}$ with $\deg(\mathfrak{g}_0)>1$ and $|\nu_\chi(\mathfrak{g}_0)|\geq1$.
Then we get
\begin{align*}
&\sum\limits_{\mathfrak{g}\in\mathfrak{G}}\bigg(|\nu_\chi(\mathfrak{g})|\cdot\deg(\mathfrak{g})-|\nu_\chi(\mathfrak{g})|\cdot\gcd(\deg(\mathfrak{g}),1)\bigg)\\
&\geq|\nu_\chi(\mathfrak{g}_0)|\cdot\deg(\mathfrak{g}_0)-|\nu_\chi(\mathfrak{g}_0)|\cdot\gcd(\deg(\mathfrak{g}_0),1)\\
&=|\nu_\chi(\mathfrak{g}_0)|\cdot\deg(\mathfrak{g}_0)-|\nu_\chi(\mathfrak{g}_0)|\\
&=|\nu_\chi(\mathfrak{g}_0)|\cdot(\deg(\mathfrak{g}_0)-1)\\
&\geq\deg(\mathfrak{g}_0)-1\\
&\geq1
\end{align*}
so by Lemma 4.3 we have $(q-1)|\chi(1)$ and since $p|(q-1)$, we have $p|\chi(1)$.
In this case Property 1.3 does not make a statement.
This shows that we can assume $\deg(\mathfrak{g})=1$ for all $\mathfrak{g}\in\mathfrak{G}$ with $|\nu_\chi(\mathfrak{g})|\geq1$ and by Lemma 2.1, we can assume that $p^2|f_p$ which then implies $p^2|(q-1)$.

Now every polynomial $\mathfrak{f}\in\mathfrak{F}$ with $\deg(\mathfrak{f})=1$ is of the form $x-z^a$ for a primitive $(q-1)$-th root of unity $z$ and some parameter $a\in\{1,\dots,q-1\}$.
We choose the conjugacy class $c$ such that
\[\nu_\chi(\mathfrak{f}_i)=\{1^1\},\]
for all $i\in\{1,\dots,n\}$ and $()$ otherwise, where $\mathfrak{f}_1,\dots,\mathfrak{f}_n\in\mathfrak{F}$ are degree $1$ polynomials with respective parameters $a_1,\dots,a_n$ which are chosen as follows
\[a_i=i\cdot p\cdot(q-1)_{p'}\text{ for }i\in\{1,\dots,n-3\},\]
\[a_{n-2}=j\cdot p\cdot(q-1)_{p'}\]
\[a:=a_{n-1}=(q-1)_{p'},\ b:=a_n=q-1-(q-1)_{p'},\]
where
\[j\in\{1,\dots,q-1\}\text{ such that }(q-1)|(a_1+\dots+a_{n-2}).\]
Note that it is possible to choose $\mathfrak{f}_1,\dots,\mathfrak{f}_n$ in this way since $p^2|(q-1)$ and therefore
\[n\cdot p\cdot(q-1)_{p'}<p^2\cdot(q-1)_{p'}\leq q-1.\]
Further, note that $a\ne b$ since $p\ne2$, so $a_1,\dots,a_n$ are pairwise different and that the elements of $c$ lie in $\SL_n(q)$ since
\[a_1+\dots+a_n\equiv0\mod(q-1),\]
so $\det(g)=1$ for all $g$ lying in the class $c$.
With this choice of $c$, we have that the elements of $c$ are $p$-elements and their order is $(q-1)_p$.

The only partition $\rho$ of $n$ for which there is a contribution to $\chi(c)$ is $\rho=\{1^n\}$, since for all other partitions of $n$ there is no mode of substitution into the class $c$.
It is then clear that for $\rho$ there is only one mode of substitution
$m_c$ into $c$ and only one mode of substitution
$m_\chi$ into $\chi$, so we get
\[\chi(c)=(-1)^{n-\sum\limits_{\mathfrak{g}\in\mathfrak{G}}|\nu_\chi(\mathfrak{g})|}\mu(m_\chi,\chi)Q(m_c,c)\cdot\]
\[\prod\limits_d\sum\limits_{\sigma\in\S_{r_d}}S_d(h_{d1}m_\chi:\xi_{d\sigma(1)})\cdots S_d(h_{dr_d}m_\chi:\xi_{d\sigma(r_d)}).\]
The terms $\mu(m_\chi,\chi)$ and $Q(m_c,c)$ are non-zero rational numbers and since we are only interested in the field $\mathbb{Q}(\chi(c))$, we can restrict attention to the term
\[\prod\limits_d\sum\limits_{\sigma\in\S_{r_d}}S_d(h_{d1}m_\chi:\xi_{d\sigma(1)})\cdots S_d(h_{dr_d}m_\chi:\xi_{d\sigma(r_d)}).\]
Now the product only has one factor, namely the one where $d=1$ and then $r_d=n$, so the term becomes
\[\sum\limits_{\sigma\in\S_n}S_1(h_{11}m_\chi:\xi_{1\sigma(1)})\cdots S_1(h_{1n}m_\chi:\xi_{1\sigma(n)}).\]
Further, each term $\theta^k(\xi)$ is of the form $z^{a_ik_j}$, for some $i,j\in\{1,\dots,n\}$ where $\{k_1\},\dots,\{k_n\}$ are the simpleces which are mapped to a non-empty partition by $\nu_\chi$.
Note that $k_1,\dots,k_n$ are not necessarily pairwise different.
The term we want to investigate can then be written as
\[S:=\sum\limits_{\sigma\in\S_n}z^{a_1k_{\sigma(1)}+\dots+a_nk_{\sigma(n)}}.\]
We can see that $S$ is the value of a sum of $f'_p$-th roots of unity for
\[f'_p:=\frac{(q-1)_p}{\gcd(q-1,k_1,\dots,k_n)_p}\]
and if $k':=\gcd(q-1,k_1,\dots,k_l)_p>1$, then
\[z^{a_1k_{\sigma(1)}+\dots+a_lk_{\sigma(n)}}=z^{k'\cdot(a_1k'_{\sigma(1)}+\dots+a_lk'_{\sigma(n)})}=y^{a_1k'_{\sigma(1)}+\dots+a_nk'_{\sigma(n)}}\]
for $k'_i:=k_i/k'$ and a primitive $f'_p$-th root of unity $y$.
In this case, we replace $z$ by $y$ and $k_i$ by $k'_i$, so we can assume that $p\nmid k_i$ for at least one $i\in\{1,\dots,n\}$.

For the next step of the proof, let $l_i$ be the number of permutations $\sigma\in\S_n$ with
\[a_1k_{\sigma(1)}+\dots+a_nk_{\sigma(n)}\equiv i\mod p\]
for $i\in\{0,\dots,p-1\}$.
Now every permutation contributes to exactly one $l_i$ and since $p|a_j$ for all $j\in\{1,\dots,n-2\}$, we have that $i$ does not depend on $\sigma(1),\dots,\sigma(n-2)$.
This implies that $l_i$ is a multiple of $(n-2)!$ and since $l_i\leq n!$, we have
\[l_i=(n-2)!\cdot b_i\]
for some $b_i\in\{0,\dots,n(n-1)\}$.
Now suppose that there is some $k\in\{1,\dots,p-1\}$ with $l_k\ne0$.
Then since $p>n(n-1)$, we have $p\nmid(n-2)!$ and $p\nmid b_k$, so $p\nmid l_k$.
It follows that if $s$ is a sum of $f'_p$-th roots of unity with $|s|=S$, then
$|s^{(f'_p)}|\ne0$ since if $z^i$ and $z^j$ are entries of the same minimal vanishing sum $t$, which has length $p$, then $i\equiv j\mod p$, so $|s^{(f'_p)}|\ne0$ is not a union of minimal vanishing sums.
By Theorem 2.5 (b) we then have $f_p=f'_p$ and Property 1.3 holds by Theorem 2.5 (c).

Now consider the case where $l_i=0$ for all $i\in\{1,\dots,p-1\}$.
This happens if and only if
\[k_1\equiv\dots\equiv k_n\mod p.\]
Here we choose the class $c$ differently by setting
\[a_i=i\cdot p\cdot(q-1)_{p'}\text{ for }i\in\{1,\dots,n-1\},\]
\[a:=a_n=(q-1)_{p'}.\]
With this choice of $c$, the elements of $c$ are again $p$-elements of order $(q-1)_p$.
The proof in this case is similar to the proof of the previous case.
Let $l'_i$ be the number of permutations $\sigma\in\S_n$ with
\[a_1k_{\sigma(1)}+\dots+a_nk_{\sigma(n)}\equiv i\mod p,\ i\in\{0,\dots,p-1\}\]
for the new $a_1,\dots,a_n$.
Then there is some $i\in\{0,\dots,p-1\}$ with
\[l'_i=n!,\ l'_j=0\]
for all $j\in\{0,\dots,p-1\},\ j\ne i$.
Note that $i\ne0$ since otherwise we would have
\[k_1\equiv\dots\equiv k_n=0\mod p,\]
which would contradict the assumption that there is some $i\in\{1,\dots,n\}$ with $p\nmid k_i$.
Then again Property 1.3. holds by Theorem 2.5 (b)+(c).
Note that in this case the elements of $c$ do not lie in $\SL_n(q)$.
\end{proof}
\end{theorem}

\begin{remark}
In the proof of Theorem 4.4 it is also possible to choose the conjugacy class $c$ like in the last part of the proof for every irreducible complex character $\chi$ with $\deg(\mathfrak{g})=1$ for all $\mathfrak{g}\in\mathfrak{G}$ with $|\nu_\chi(\mathfrak{g})|\geq1$.
With this choice the $p$-element $g$ never lies in $\SL_n(q)$, but since 
$l_i=n!$ for some $i\in\{1,\dots,p-1\}$, we have $p\nmid l_i$ for all primes $p$ with $p>n$ and then Property 1.3 holds by Theorem 2.5 (b)+(c).
This shows that Property 1.3 holds for the finite general linear groups for primes $p$ with $p|(q-1)$ and $p>n$.
\end{remark}

Another question to ask is the following.
If $H$ is a group which has Property 1.3 and $N\unlhd H$ is a normal subgroup of $H$, under which assumptions does $N$ also have Property 1.3?
The following theorem deals with this situation.

\begin{theorem}
Let $G$ be a finite group and let $N\unlhd G$.
Further, let $\chi\in\Irr(G)$ and $\psi\in\Irr(N)$ such that $\psi$ is a constituent of $\chi\downarrow_N^G$ and let $p$ be an odd prime with $p\nmid[G:\mathcal{I}_G(\psi)]$ and $p\nmid e$, where
\[\mathcal{I}_G(\psi):=\{g\in G\ |\ {^g}\psi=\psi\}\]
is the inertia group of $\psi$ in $G$ and $e$ is the ramification index of $\chi$ in $N$ (see \cite[Section 6]{Isa}).
Assume that Property 1.3 holds for $G,\ \chi$ and $p$ and that the $p$-element $g$ for which $p\nmid[\mathbb{Q}_{(f_{\chi})_p}:\mathbb{Q}(\chi(g))]$ can be chosen such that $g\in N$.
Then Property 1.3 also holds for $N,\ \psi$ and $p$.
\begin{proof}
By \cite[Theorem 6.2]{Isa} we have
\[\chi\downarrow_N^G=e\sum\limits_{g\in G/\mathcal{I}_G(\psi)}{^g}\psi\]
for some $e\in\mathbb{Z}_{\geq1}$, where
\[{^g}\psi:N\rightarrow\mathbb{C},\ x\mapsto\psi(g^{-1}xg),\]
is the conjugate of $\psi$ by the element $g$.
If $(f_{\psi})_p\leq p$, then Property 1.3 holds for $\psi$ by Lemma 2.1, so assume that $(f_{\psi})_p\geq p^2$.
Further, we can assume that $p\nmid\psi(1)$.
Since all conjugates of $\psi$ have the same degree, we have
\[\chi(1)=e\cdot[G:\mathcal{I}_G(\psi)]\cdot\psi(1)\]
and the assumptions imply $p\nmid\chi(1)$.

Now if $(f_{\chi})_p>1$, then by \cite[Lemma 4.2 (ii)]{NaTi} we have $(f_{\psi})_p\leq (f_{\chi})_p$ since for this claim of the lemma the assumption $p\nmid[G:\mathcal{I}_G(\psi)]$ suffices.
On the other hand, if $(f_{\chi})_p=1$, then $\chi$ is $\sigma$-invariant for $\sigma\in\Gal(\mathbb{Q}^\ab/\mathbb{Q})$, which fixes roots of unity of order not divisible by $p$ and maps any $p$-power root of unity $z$ to $z^{p+1}$.
It follows that $(f_{\psi})_p\leq p$ by \cite[Lemma 4.1 (iii) + Lemma 4.2 (i)]{NaTi} which is a contradiction to $(f_{\psi})_p\geq p^2$.

Then since Property 1.3 holds for $G,\ \chi$ and $p$, there is a $p$-element $h\in G$ for which we have $p\nmid[\mathbb{Q}_{(f_{\chi})_p}:\mathbb{Q}(\chi(h))]$.
Further, we can assume that $h\in N$ by assumption.
Now let $s,\ t_g$ for $g\in G$ be sums of $(f_{\chi})_p$-th roots of unity with $|s|=\chi(h),\ |t_g|={^g}\psi(h)$.
Then
\[|s|=e\cdot\sum\limits_{g\in G/\mathcal{I}_G(\psi)}|t_g|\text{ and }|s^{((f_{\chi})_p)}|=e\cdot\sum\limits_{g\in G/\mathcal{I}_G(\psi)}|t_g^{((f_{\chi})_p)}|\]
and since $p\nmid[\mathbb{Q}_{(f_{\chi})_p}:\mathbb{Q}(s)]$, we have $|s^{((f_{\chi})_p)}|\ne0$ by Theorem 2.5 (b)+(c).
But then there exists some $g_0\in G$ with $|t_{g_0}^{((f_{\chi})_p)}|\ne0$ and since
\[{^{g_0}}\psi(h)=\psi(g_0^{-1}hg_0),\]
we have that Property 1.3 holds for $N$ and $\psi$ by choosing $g_0^{-1}hg_0$ as the $p$-element and applying Theorem 2.5 (b)+(c).
\end{proof}
\end{theorem}

We can apply Theorem 4.6 to the case where $G$ is the general linear group in dimension $n$ and the normal subgroup $N$ is the associated special linear group.
This is done in the following theorem.

\begin{theorem}
Let $n\in\mathbb{Z}_{\geq1}$, let $q$ be a prime power with $(n,q)\notin\{(2,2),(2,3)\}$ and let $p$ be a prime with $p|(q-1)$ and $p>n(n-1)$.
Then Property 1.3 holds for all quasi-simple groups $G$ with $S:=G/\Z(G)=\PSL_n(q)$ for the prime $p$.
\begin{proof}
By \cite[Section 24.2]{MaTe} and by the proof of Theorem 3.1 it is enough to show that Property 1.3 holds for $\SL_n(q)$ for the prime $p$.

We need to verify the assumptions of Theorem 4.6.
Set $G:=\GL_n(q),\ N:=\SL_n(q)$ and let $\psi\in\Irr(N)$.
Further, let $\chi\in\Irr(G)$ be an irreducible character with conductor $f$ such that $\psi$ is a summand of $\chi\downarrow_N^G$ and let $k$ be the number of irreducible summands of $\chi\downarrow_N^G$.

We can assume that $\deg(\mathfrak{g})=1$ for all $\mathfrak{g}\in\mathfrak{G}$ with $|\nu_\chi(\mathfrak{g})|\geq1$.
Indeed, if $\mathfrak{g}_1,\dots,\mathfrak{g}_l$ are the simpleces of degree greater than $1$ which are mapped to a non-empty partition by $\nu_\chi$, then
\begin{align*}
&\sum\limits_{\mathfrak{g}\in\mathfrak{G}}\bigg(|\nu_\chi(\mathfrak{g})|\cdot\deg(\mathfrak{g})-\Big\lfloor|\nu_\chi(\mathfrak{g})|\cdot\gcd(\deg(\mathfrak{g}),1)\Big\rfloor_1\bigg)\\
&=\sum\limits_{\mathfrak{g}\in\mathfrak{G}}\bigg(|\nu_\chi(\mathfrak{g})|\cdot\deg(\mathfrak{g})-|\nu_\chi(\mathfrak{g})|\cdot\gcd(\deg(\mathfrak{g}),1)\bigg)\\
&=\sum\limits_{i=1}^l\bigg(|\nu_\chi(\mathfrak{g}_i)|\cdot\deg(\mathfrak{g}_i)-|\nu_\chi(\mathfrak{g}_i)|\cdot\gcd(\deg(\mathfrak{g}_i),1)\bigg)
\end{align*}
which is greater than $1$ unless $l=1,\ \deg(\mathfrak{g}_1)=2$ and $\nu_\chi(\mathfrak{g}_1)=\{1^1\}$.
However, in this case it follows from \cite{KlTi}, that $k=1$ and then it is clear that Property 1.3 holds for $\psi$.
So we have that $(q-1)^2|\chi(1)$ by Lemma 4.3 and since $[G:N]=q-1$, it follows that $(q-1)|\psi(1)$.
This implies that $p|\psi(1)$ and then Property 1.3 does not make a statement.

Now it follows from \cite[Theorem 1.1]{KlTi} that if $\rho$ is a partition in the image of $\nu_\chi$, then the number of simpleces which are mapped to $\rho$ by $\nu_\chi$ is divisible by $k$ and since
\[\sum\limits_{\mathfrak{g}\in\mathfrak{G}}|\nu_\chi(\mathfrak{g})|\cdot\deg(\mathfrak{g})=n,\]
we have $k|n$ and therefore $k|\gcd(n,q-1)$.
The assumption $p>n(n-1)$, then implies $p\nmid k$.
Now since $G/N$ is cyclic, we have $k=[G:\mathcal{I_G(\psi)}]$ and the ramification index of $\chi$ in $N$ is equal to $1$ by \cite[Section 3]{KlTi}.
This shows two assumptions of Theorem 4.6.

Now by Theorem 4.4 we have that Property 1.3 holds for $G,\ \chi$ and $p$.
It is left to show that we can choose the $p$-element $g$ for which $p\nmid[\mathbb{Q}_{f_p}:\mathbb{Q}(\chi(g))]$ such that $g\in N$.

As in the proof of Theorem 4.4, we can assume that there is some $i\in\{1,\dots,n\}$ with $p\nmid k_i$ by replacing $k_1,\dots,k_n$ with $k'_1,\dots,k'_n$, where $k'_i:=k_i/k'$ for all $i\in\{1,\dots,n\}$ and $k':=\gcd(q-1,k_1,\dots,k_n)_p$.

By the proof of Theorem 4.4 we can choose the $p$-element $g$ for which we have $p\nmid[\mathbb{Q}_{f_p}:\mathbb{Q}(\chi(g))]$ such that $g\in N$, unless
\[k_1\equiv\dots\equiv k_n\mod p,\]
where $k_1,\dots,k_n$ are the parameters of $\chi$.
In this case set
\[m:=\max\{i\in\mathbb{Z}_{\geq1}\ |\ k_1\equiv\dots\equiv k_n\mod p^i\}\in\mathbb{Z}_{\geq1}\cup\{\infty\}.\]
Now for an arbitrary $h\in N$ in a conjugacy class $c_h$ with $\deg(\mathfrak{g})=1$ for all $\mathfrak{g}\in\mathfrak{G}$ with $|\nu_{c_h}(\mathfrak{g})|=1$ and parameters $b_1,\dots,b_n$ we have
\[b_1+\dots+b_n\equiv0\mod(q-1)\]
since $\det(h)=1$ and therefore
\[p^i|(b_1k_{\sigma(1)}+\dots+b_nk_{\sigma(n)})\]
for all $i\in\{1,\dots,m\}$ with $p^i\leq(q-1)_p$ and all permutations $\sigma\in\S_n$.
It follows that if $p^{m+1}\geq(q-1)_p/k'$, then $f_p\leq p$ and therefore Property 1.3 holds by Lemma 2.1.

So assume $p^{m+1}<(q-1)_p/k'$.
Then $\chi(h)$ is a sum of roots of unity of order not divisible by
\[\frac{(q-1)_p}{k'\cdot p^{m-1}},\]
so we have
\[f_p\leq f'_p:=\frac{(q-1)_p}{k'\cdot p^m}.\]
We choose the conjugacy class $c$ such that the parameters $a_1,\dots,a_n$ of $c$ are as follows
\[a_i=i\cdot p^{m+1}\cdot(q-1)_{p'}\text{ for }i\in\{1,\dots,n-3\},\]
\[a_{n-2}=j\cdot p^{m+1}\cdot(q-1)_{p'},\]
\[a_{n-1}=(q-1)_{p'},\ a_n=q-1-(q-1)_{p'},\]
where
\[j\in\{1,\dots,q-1\}\text{ such that }(q-1)|(a_1+\dots+a_{n-2}).\]
Now let $l_i$ be the number of permutations $\sigma\in\S_n$ with
\[a_1k_{\sigma(1)}+\dots+a_nk_{\sigma(n)}\equiv i\cdot p^m\mod p^{m+1}\]
for $i\in\{0,\dots,p-1\}$.
Then by the definition of $m$ there is some $i_0\in\{1,\dots,p-1\}$ with $l_{i_0}\ne0$ and as in the proof of Theorem 4.4 we have that
\[l_{i_0}=(n-2)!\cdot b_i\]
for some $b_i\in\{1,\dots,n(n-1)\}$, so $p\nmid l_{i_0}$ since $p>n(n-1)$.
It follows that $f_p=f'_p$ by Theorem 2.5 (b) and then Property 1.3 holds by Theorem 2.5 (c).
\end{proof}
\end{theorem}

\section{Tables of $p$-elements}

\centering
\begin{tabular}{|c|c|c|c|c|}
\multicolumn{5}{c}{Table 1 \ $p$-elements for $\SL_2(q),\ q$ odd (for the character}\\
\multicolumn{5}{c}{table in \cite{CHEVIE} which was computed in \cite{Sch}, see also \cite[Table 5.4]{Bon})}\\
\hline
character $\chi$&polynomial $Q$&$p$-power $f'_p$&class of $g$&$\Tilde{\lambda}(\chi(g),f'_p)$\\
\hline
\hline
$\chi_4(k)$&$q-1$&$\frac{Q_p}{\gcd(Q,k)_p}$&$C_2(Q_{p'})$&$2$\\
\hline
$\chi_5(k)$&$q+1$&$\frac{Q_p}{\gcd(Q,k)_p}$&$C_3(Q_{p'})$&$2$\\
\hline
\end{tabular}\\
$\ $\\
$\ $\\
\begin{tabular}{|c|c|c|c|c|}
\multicolumn{5}{c}{Table 2 \ $p$-elements for $\SL_2(q),\ q$ even}\\
\multicolumn{5}{c}{(for the character table in \cite{CHEVIE} which was computed in \cite{Sch})}\\
\hline
character $\chi$&polynomial $Q$&$p$-power $f'_p$&class of $g$&$\Tilde{\lambda}(\chi(g),f'_p)$\\
\hline
\hline
$\chi_2(n)$&$q-1$&$\frac{Q_p}{\gcd(Q,n)_p}$&$C_2(Q_{p'})$&$2$\\
\hline
$\chi_3(n)$&$q+1$&$\frac{Q_p}{\gcd(Q,n)_p}$&$C_3(Q_{p'})$&$2$\\
\hline
\end{tabular}\\
$\ $\\
$\ $\\
\begin{tabular}{|c|c|c|c|c|}
\multicolumn{5}{c}{Table 3 \ $p$-elements for $\SL_3(q)$ (for \cite[Table 1b]{SiFr})}\\
\hline
character $\chi$&polynomial $Q$&$p$-power $f'_p$&class of $g$&$\Tilde{\lambda}(\chi(g),f'_p)$\\
\hline
\hline
$\chi_{rt}^{(u)}$&$q+1$&$\frac{Q_p}{\gcd(Q,u)_p}$&$C_7^{((q-1)\cdot Q_{p'})}$&$2$\\
\hline
$\chi_t^{(u)}$&$q-1$&$\frac{Q_p}{\gcd(Q,u)_p}$&$C_6^{(Q_{p'},Q-Q_{p'},Q)}$&$2$\\
\hline
$\chi_{qt}^{(u)}$&$q-1$&$\frac{Q_p}{\gcd(Q,u)_p}$&$C_6^{(Q_{p'},Q-Q_{p'},Q)}$&$2$\\
\hline
$\chi_{st}^{(u,v,w)}$&$q-1$&$\frac{Q_p}{\gcd(Q,\lcm(u,v,w))_p}$&$C_4^{(Q_{p'})}$&$\{1,2,3\}$\\
\hline
$\chi_{r^2s}^{(u)}$&$q^2+q+1$&$\frac{Q_p}{\gcd(Q,u)_p}$&$C_8^{(Q_{p'})}$&$3$\\
\hline
\end{tabular}\\
$\ $\\
$\ $\\
\begin{tabular}{|c|c|c|c|c|}
\multicolumn{5}{c}{Table 4 \ $p$-elements for $\SU_3(q)$ (for \cite[Table 1b]{SiFr})}\\
\hline
character $\chi$&polynomial $Q$&$p$-power $f'_p$&class of $g$&$\Tilde{\lambda}(\chi(g),f'_p)$\\
\hline
\hline
$\chi_{rt}^{(u)}$&$q-1$&$\frac{Q_p}{\gcd(Q,u)_p}$&$C_7^{((q+1)\cdot Q_{p'})}$&$2$\\
\hline
$\chi_t^{(u)}$&$q+1$&$\frac{Q_p}{\gcd(Q,u)_p}$&$C_6^{(Q_{p'},Q-Q_{p'},Q)}$&$2$\\
\hline
$\chi_{qt}^{(u)}$&$q+1$&$\frac{Q_p}{\gcd(Q,u)_p}$&$C_6^{(Q_{p'},Q-Q_{p'},Q)}$&$2$\\
\hline
$\chi_{st}^{(u,v,w)}$&$q+1$&$\frac{Q_p}{\gcd(Q,\lcm(u,v,w))_p}$&$C_4^{(Q_{p'})}$&$\{1,2,3\}$\\
\hline
$\chi_{r^2s}^{(u)}$&$q^2-q+1$&$\frac{Q_p}{\gcd(Q,u)_p}$&$C_8^{(Q_{p'})}$&$3$\\
\hline
\end{tabular}\\
$\ $\\
$\ $\\
\begin{tabular}{|c|c|c|c|c|}
\multicolumn{5}{c}{Table 5 \ $p$-elements for $\Sp_4(q),\ q$ odd}\\
\multicolumn{5}{c}{(for the character table in \cite{Sri})}\\
\hline
character $\chi$&polynomial $Q$&$p$-power $f'_p$&class of $g$&$\Tilde{\lambda}(\chi(g),f'_p)$\\
\hline
\hline
$\chi_1(j)$&$q^2+1$&$\frac{Q_p}{\gcd(Q,j)_p}$&$B_1(Q_{p'})$&$4$\\
\hline
$\chi_3(k,l)$&$q-1$&$\frac{Q_p}{\gcd(Q,\lcm(k,l))_p}$&$C_3(Q_{p'})$&$\{2,4\}$\\
\hline
$\chi_4(k,l)$&$q+1$&$\frac{Q_p}{\gcd(Q,\lcm(k,l))_p}$&$C_1(Q_{p'})$&$\{2,4\}$\\
\hline
$\chi_6(k)$&$q+1$&$\frac{Q_p}{\gcd(Q,k)_p}$&$C_1(Q_{p'})$&$2$\\
\hline
$\chi_7(k)$&$q+1$&$\frac{Q_p}{\gcd(Q,k)_p}$&$C_1(Q_{p'})$&$2$\\
\hline
$\chi_8(k)$&$q-1$&$\frac{Q_p}{\gcd(Q,k)_p}$&$C_3(Q_{p'})$&$2$\\
\hline
$\chi_9(k)$&$q-1$&$\frac{Q_p}{\gcd(Q,k)_p}$&$C_3(Q_{p'})$&$2$\\
\hline
$\xi_1(k)$&$q+1$&$\frac{Q_p}{\gcd(Q,k)_p}$&$B_6(Q_{p'})$&$2$\\
\hline
$\xi'_1(k)$&$q+1$&$\frac{Q_p}{\gcd(Q,k)_p}$&$B_6(Q_{p'})$&$2$\\
\hline
$\xi_3(k)$&$q-1$&$\frac{Q_p}{\gcd(Q,k)_p}$&$B_8(Q_{p'})$&$2$\\
\hline
$\xi'_3(k)$&$q-1$&$\frac{Q_p}{\gcd(Q,k)_p}$&$B_8(Q_{p'})$&$2$\\
\hline
$\xi'_{21}(k)$&$q+1$&$\frac{Q_p}{\gcd(Q,k)_p}$&$B_6(Q_{p'})$&$2$\\
\hline
$\xi_{41}(k)$&$q-1$&$\frac{Q_p}{\gcd(Q,k)_p}$&$B_8(Q_{p'})$&$2$\\
\hline
\end{tabular}\\
$\ $\\
$\ $\\
\begin{tabular}{|c|c|c|c|c|}
\multicolumn{5}{c}{Table 6 \ $p$-elements for $\Sp_4(q),\ q$ even (for \cite[Table IV-2]{Eno2})}\\
\hline
character $\chi$&polynomial $Q$&$p$-power $f'_p$&class of $g$&$\Tilde{\lambda}(\chi(g),f'_p)$\\
\hline
\hline
$\chi_1(k,l)$&$q-1$&$\frac{Q_p}{\gcd(Q,\lcm(k,l))_p}$&$C_1(Q_{p'})$&$\{2,4\}$\\
\hline
$\chi_4(k,l)$&$q+1$&$\frac{Q_p}{\gcd(Q,\lcm(k,l))_p}$&$C_3(Q_{p'})$&$\{2,4\}$\\
\hline
$\chi_5(k)$&$q^2+1$&$\frac{Q_p}{\gcd(Q,k)_p}$&$B_5(Q_{p'})$&$4$\\
\hline
$\chi_6(k)$&$q-1$&$\frac{Q_p}{\gcd(Q,k)_p}$&$C_1(Q_{p'})$&$2$\\
\hline
$\chi_7(k)$&$q-1$&$\frac{Q_p}{\gcd(Q,k)_p}$&$C_2(Q_{p'})$&$2$\\
\hline
$\chi_8(k)$&$q+1$&$\frac{Q_p}{\gcd(Q,k)_p}$&$C_3(Q_{p'})$&$2$\\
\hline
$\chi_9(k)$&$q+1$&$\frac{Q_p}{\gcd(Q,k)_p}$&$C_4(Q_{p'})$&$2$\\
\hline
$\chi_{10}(k)$&$q-1$&$\frac{Q_p}{\gcd(Q,k)_p}$&$C_1(Q_{p'})$&$2$\\
\hline
$\chi_{11}(k)$&$q-1$&$\frac{Q_p}{\gcd(Q,k)_p}$&$C_2(Q_{p'})$&$2$\\
\hline
$\chi_{12}(k)$&$q+1$&$\frac{Q_p}{\gcd(Q,k)_p}$&$C_3(Q_{p'})$&$2$\\
\hline
$\chi_{13}(k)$&$q+1$&$\frac{Q_p}{\gcd(Q,k)_p}$&$C_4(Q_{p'})$&$2$\\
\hline
\end{tabular}\\
$\ $\\
$\ $\\
\begin{tabular}{|c|c|c|c|c|}
\multicolumn{5}{c}{Table 7 \ $p$-elements for $\Sp_6(q),\ q$ even}\\
\multicolumn{5}{c}{(for the character table in \cite{CHEVIE} which was computed in \cite{Loo})}\\
\hline
character $\chi$&polynomial $Q$&$p$-power $f'_p$&class of $g$&$\Tilde{\lambda}(\chi(g),f'_p)$\\
\hline
\hline
$\chi_{6,1}(k)$&$q-1$&$\frac{Q_p}{\gcd(Q,k)_p}$&$C_3(Q_{p'})$&$2$\\
\hline
$\chi_{6,2}(k)$&$q-1$&$\frac{Q_p}{\gcd(Q,k)_p}$&$C_3(Q_{p'})$&$2$\\
\hline
$\chi_{6,3}(k)$&$q-1$&$\frac{Q_p}{\gcd(Q,k)_p}$&$C_3(Q_{p'})$&$2$\\
\hline
$\chi_{6,4}(k)$&$q-1$&$\frac{Q_p}{\gcd(Q,k)_p}$&$C_3(Q_{p'})$&$2$\\
\hline
$\chi_{6,6}(k)$&$q-1$&$\frac{Q_p}{\gcd(Q,k)_p}$&$C_3(Q_{p'})$&$2$\\
\hline
$\chi_{7,1}(k)$&$q+1$&$\frac{Q_p}{\gcd(Q,k)_p}$&$C_4(Q_{p'})$&$2$\\
\hline
$\chi_{7,3}(k)$&$q+1$&$\frac{Q_p}{\gcd(Q,k)_p}$&$C_4(Q_{p'})$&$2$\\
\hline
$\chi_{7,4}(k)$&$q+1$&$\frac{Q_p}{\gcd(Q,k)_p}$&$C_4(Q_{p'})$&$2$\\
\hline
$\chi_{7,5}(k)$&$q+1$&$\frac{Q_p}{\gcd(Q,k)_p}$&$C_4(Q_{p'})$&$2$\\
\hline
$\chi_{7,6}(k)$&$q+1$&$\frac{Q_p}{\gcd(Q,k)_p}$&$C_4(Q_{p'})$&$2$\\
\hline
$\chi_{8,1}(k)$&$q-1$&$\frac{Q_p}{\gcd(Q,k)_p}$&$C_3(Q_{p'})$&$2$\\
\hline
$\chi_{8,2}(k)$&$q-1$&$\frac{Q_p}{\gcd(Q,k)_p}$&$C_3(Q_{p'})$&$2$\\
\hline
$\chi_{8,3}(k)$&$q-1$&$\frac{Q_p}{\gcd(Q,k)_p}$&$C_3(Q_{p'})$&$2$\\
\hline
$\chi_{9,1}(k)$&$q+1$&$\frac{Q_p}{\gcd(Q,k)_p}$&$C_4(Q_{p'})$&$2$\\
\hline
$\chi_{9,2}(k)$&$q+1$&$\frac{Q_p}{\gcd(Q,k)_p}$&$C_4(Q_{p'})$&$2$\\
\hline
$\chi_{9,3}(k)$&$q+1$&$\frac{Q_p}{\gcd(Q,k)_p}$&$C_4(Q_{p'})$&$2$\\
\hline
$\chi_{11,1}(k)$&$q-1$&$\frac{Q_p}{\gcd(Q,k)_p}$&$C_3(Q_{p'})$&$2$\\
\hline
$\chi_{11,2}(k)$&$q-1$&$\frac{Q_p}{\gcd(Q,k)_p}$&$C_3(Q_{p'})$&$2$\\
\hline
$\chi_{11,3}(k)$&$q-1$&$\frac{Q_p}{\gcd(Q,k)_p}$&$C_3(Q_{p'})$&$2$\\
\hline
$\chi_{11,4}(k)$&$q-1$&$\frac{Q_p}{\gcd(Q,k)_p}$&$C_3(Q_{p'})$&$2$\\
\hline
$\chi_{13,1}(k)$&$q+1$&$\frac{Q_p}{\gcd(Q,k)_p}$&$C_4(Q_{p'})$&$2$\\
\hline
$\chi_{13,2}(k)$&$q+1$&$\frac{Q_p}{\gcd(Q,k)_p}$&$C_4(Q_{p'})$&$2$\\
\hline
$\chi_{13,3}(k)$&$q+1$&$\frac{Q_p}{\gcd(Q,k)_p}$&$C_4(Q_{p'})$&$2$\\
\hline
$\chi_{13,4}(k)$&$q+1$&$\frac{Q_p}{\gcd(Q,k)_p}$&$C_4(Q_{p'})$&$2$\\
\hline
$\chi_{16,1}(k,l)$&$q-1$&$\frac{Q_p}{\gcd(Q,\lcm(k,l))_p}$&$C_3(Q_{p'})$&$\{2,4\}$\\
\hline
$\chi_{16,2}(k,l)$&$q-1$&$\frac{Q_p}{\gcd(Q,\lcm(k,l))_p}$&$C_3(Q_{p'})$&$\{2,4\}$\\
\hline
$\chi_{17,1}(k,l)$&$q-1$&$\frac{Q_p}{\gcd(Q,\lcm(k,l))_p}$&$C_3(Q_{p'})$&$\{2,4\}$\\
\hline
$\chi_{17,2}(k,l)$&$q-1$&$\frac{Q_p}{\gcd(Q,\lcm(k,l))_p}$&$C_3(Q_{p'})$&$\{2,4\}$\\
\hline
$\chi_{22,1}(k,l)$&$q+1$&$\frac{Q_p}{\gcd(Q,\lcm(k,l))_p}$&$C_4(Q_{p'})$&$\{2,4\}$\\
\hline
$\chi_{22,2}(k,l)$&$q+1$&$\frac{Q_p}{\gcd(Q,\lcm(k,l))_p}$&$C_4(Q_{p'})$&$\{2,4\}$\\
\hline
$\chi_{23,1}(k,l)$&$q+1$&$\frac{Q_p}{\gcd(Q,\lcm(k,l))_p}$&$C_4(Q_{p'})$&$\{2,4\}$\\
\hline
$\chi_{23,2}(k,l)$&$q+1$&$\frac{Q_p}{\gcd(Q,\lcm(k,l))_p}$&$C_4(Q_{p'})$&$\{2,4\}$\\
\hline
$\chi_{24,1}(k)$&$q^2+1$&$\frac{Q_p}{\gcd(Q,k)_p}$&$C_{20}(Q_{p'})$&$4$\\
\hline
$\chi_{24,2}(k)$&$q^2+1$&$\frac{Q_p}{\gcd(Q,k)_p}$&$C_{20}(Q_{p'})$&$4$\\
\hline
$\chi_{25,1}(k,l,m)$&$q-1$&$\frac{Q_p}{\gcd(Q,\lcm(k,l,m))_p}$&$C_3(Q_{p'})$&$\{2,4,6\}$\\
\hline
$\chi_{30,1}(k,l)$&$q^2+1$&$\frac{Q_p}{\gcd(Q,l)_p}$&$C_{20}(Q_{p'})$&$4$\\
\hline
$\chi_{31,1}(k)$&$q^2+q+1$&$\frac{Q_p}{\gcd(Q,k)_p}$&$C_{28}((q-1)\cdot Q_{p'})$&$6$\\
\hline
$\chi_{32,1}(k,l,m)$&$q+1$&$\frac{Q_p}{\gcd(Q,\lcm(k,l,m))_p}$&$C_4(Q_{p'})$&$\{2,4,6\}$\\
\hline
$\chi_{33,1}(k,l)$&$q^2+1$&$\frac{Q_p}{\gcd(Q,k)_p}$&$C_{20}(Q_{p'})$&$4$\\
\hline
$\chi_{34,1}(k)$&$q^2-q+1$&$\frac{Q_p}{\gcd(Q,k)_p}$&$C_{31}((q+1)\cdot Q_{p'})$&$6$\\
\hline
\end{tabular}\\
$\ $\\
$\ $\\
\begin{tabular}{|c|c|c|c|c|}
\multicolumn{5}{c}{Table 8 \ $p$-elements for $G_2(q),\ q$ even}\\
\multicolumn{5}{c}{(for \cite[Table IV-2 and Table IV-3]{EnYa})}\\
\hline
character $\chi$&polynomial $Q$&$p$-power $f'_p$&class of $g$&$\Tilde{\lambda}(\chi(g),f'_p)$\\
\hline
\hline
$\chi_1(k)$&$q-1$&$\frac{Q_p}{\gcd(Q,k)_p}$&$C_{11}(Q_{p'})$&$2$\\
\hline
$\chi'_1(k)$&$q+1$&$\frac{Q_p}{\gcd(Q,k)_p}$&$D_{11}(Q_{p'})$&$2$\\
\hline
$\chi_2(k)$&$q-1$&$\frac{Q_p}{\gcd(Q,k)_p}$&$C_{11}(Q_{p'})$&$2$\\
\hline
$\chi'_2(k)$&$q+1$&$\frac{Q_p}{\gcd(Q,k)_p}$&$D_{11}(Q_{p'})$&$2$\\
\hline
$\chi_3(k)$&$q-1$&$\frac{Q_p}{\gcd(Q,k)_p}$&$C_{21}(Q_{p'})$&$2$\\
\hline
$\chi'_3(k)$&$q+1$&$\frac{Q_p}{\gcd(Q,k)_p}$&$D_{21}(Q_{p'})$&$2$\\
\hline
$\chi_4(k)$&$q-1$&$\frac{Q_p}{\gcd(Q,k)_p}$&$C_{21}(Q_{p'})$&$2$\\
\hline
$\chi'_4(k)$&$q+1$&$\frac{Q_p}{\gcd(Q,k)_p}$&$D_{21}(Q_{p'})$&$2$\\
\hline
$\chi_5(k,l)$&$q-1$&$\frac{Q_p}{\gcd(Q,\lcm(k,l))_p}$&$C_{21}(Q_{p'})$&$\{4,6\}$\\
\hline
$\chi'_5(k,l)$&$q+1$&$\frac{Q_p}{\gcd(Q,\lcm(k,l))_p}$&$D_{21}(Q_{p'})$&$\{4,6\}$\\
\hline
$\chi_7(k)$&$q^2+q+1$&$\frac{Q_p}{\gcd(Q,k)_p}$&$E_3(Q_{p'})$&$6$\\
\hline
$\chi'_7(k)$&$q^2-q+1$&$\frac{Q_p}{\gcd(Q,k)_p}$&$E_4(Q_{p'})$&$6$\\
\hline
\end{tabular}\\
$\ $\\
$\ $\\
\begin{tabular}{|c|c|c|c|c|}
\multicolumn{5}{c}{Table 9 \ $p$-elements for $G_2(q),\ q=3^n$ (for \cite[Table VII-2]{Eno})}\\
\hline
character $\chi$&polynomial $Q$&$p$-power $f'_p$&class of $g$&$\Tilde{\lambda}(\chi(g),f'_p)$\\
\hline
\hline
$\chi_1(k)$&$q-1$&$\frac{Q_p}{\gcd(Q,k)_p}$&$C_{11}(Q_{p'})$&$2$\\
\hline
$\chi_2(k)$&$q-1$&$\frac{Q_p}{\gcd(Q,k)_p}$&$C_{11}(Q_{p'})$&$2$\\
\hline
$\chi_3(k)$&$q-1$&$\frac{Q_p}{\gcd(Q,k)_p}$&$C_{21}(Q_{p'})$&$2$\\
\hline
$\chi_4(k)$&$q-1$&$\frac{Q_p}{\gcd(Q,k)_p}$&$C_{21}(Q_{p'})$&$2$\\
\hline
$\chi_5(k)$&$q+1$&$\frac{Q_p}{\gcd(Q,k)_p}$&$D_{11}(Q_{p'})$&$2$\\
\hline
$\chi_6(k)$&$q+1$&$\frac{Q_p}{\gcd(Q,k)_p}$&$D_{11}(Q_{p'})$&$2$\\
\hline
$\chi_7(k)$&$q+1$&$\frac{Q_p}{\gcd(Q,k)_p}$&$D_{21}(Q_{p'})$&$2$\\
\hline
$\chi_8(k)$&$q+1$&$\frac{Q_p}{\gcd(Q,k)_p}$&$D_{21}(Q_{p'})$&$2$\\
\hline
$\chi_9(k,l)$&$q-1$&$\frac{Q_p}{\gcd(Q,\lcm(k,l))_p}$&$C_{21}(Q_{p'})$&$\{4,6\}$\\
\hline
$\chi_{12}(k,l)$&$q+1$&$\frac{Q_p}{\gcd(Q,\lcm(k,l))_p}$&$D_{21}(Q_{p'})$&$\{4,6\}$\\
\hline
$\chi_{13}(k)$&$q^2+q+1$&$\frac{Q_p}{\gcd(Q,k)_p}$&$E_5(Q_{p'})$&$6$\\
\hline
$\chi_{14}(k)$&$q^2-q+1$&$\frac{Q_p}{\gcd(Q,k)_p}$&$E_6(Q_{p'})$&$6$\\
\hline
\end{tabular}\\
$\ $\\
$\ $\\
\begin{tabular}{|c|c|c|c|c|}
\multicolumn{5}{c}{Table 10 \ $p$-elements for $G_2(q),\ q=r^n$ with $r>3$ (for}\\
\multicolumn{5}{c}{\cite[Table B.2 and Table B.3]{Hiß}, which were originally computed in \cite{ChRe})}\\
\hline
character $\chi$&polynomial $Q$&$p$-power $f'_p$&class of $g$&$\Tilde{\lambda}(\chi(g),f'_p)$\\
\hline
\hline
$X_{1a}(k)$&$q-1$&$\frac{Q_p}{\gcd(Q,k)_p}$&$h_{1b}(Q_{p'})$&$2$\\
\hline
$X'_{1a}(k)$&$q-1$&$\frac{Q_p}{\gcd(Q,k)_p}$&$h_{1b}(Q_{p'})$&$2$\\
\hline
$X_{1b}(k)$&$q-1$&$\frac{Q_p}{\gcd(Q,k)_p}$&$h_{1a}(Q_{p'})$&$2$\\
\hline
$X'_{1b}(k)$&$q-1$&$\frac{Q_p}{\gcd(Q,k)_p}$&$h_{1a}(Q_{p'})$&$2$\\
\hline
$X_{2a}(k)$&$q+1$&$\frac{Q_p}{\gcd(Q,k)_p}$&$h_{2b}(Q_{p'})$&$2$\\
\hline
$X'_{2a}(k)$&$q+1$&$\frac{Q_p}{\gcd(Q,k)_p}$&$h_{2b}(Q_{p'})$&$2$\\
\hline
$X_{2b}(k)$&$q+1$&$\frac{Q_p}{\gcd(Q,k)_p}$&$h_{2a}(Q_{p'})$&$2$\\
\hline
$X'_{2b}(k)$&$q+1$&$\frac{Q_p}{\gcd(Q,k)_p}$&$h_{2a}(Q_{p'})$&$2$\\
\hline
$X_1(k,l)$&$q-1$&$\frac{Q_p}{\gcd(Q,\lcm(k,l))_p}$&$h_{1a}(Q_{p'})$&$\{4,6\}$\\
\hline
$X_2(k,l)$&$q+1$&$\frac{Q_p}{\gcd(Q,\lcm(k,l))_p}$&$h_{2b}(Q_{p'})$&$\{4,6\}$\\
\hline
$X_3(k)$&$q^2+q+1$&$\frac{Q_p}{\gcd(Q,k)_p}$&$h_3(Q_{p'})$&$6$\\
\hline
$X_6(k)$&$q^2-q+1$&$\frac{Q_p}{\gcd(Q,k)_p}$&$h_6(Q_{p'})$&$6$\\
\hline
\end{tabular}\\
$\ $\\
$\ $\\
\begin{tabular}{|c|c|c|c|c|}
\multicolumn{5}{c}{Table 11 \ $p$-elements for ${^3}D_4(q)$}\\
\multicolumn{5}{c}{(for the character tables in \cite{CHEVIE} which were computed in \cite{DeMi})}\\
\hline
character $\chi$&polynomial $Q$&$p$-power $f'_p$&class of $g$&$\Tilde{\lambda}(\chi(g),f'_p)$\\
\hline
\hline
$\chi_{3,0}(k)$&$q-1$&$\frac{Q_p}{\gcd(Q,k)_p}$&$C_6(q^3-1,Q_{p'})$&$6$\\
\hline
$\chi_{3,1}(k)$&$q-1$&$\frac{Q_p}{\gcd(Q,k)_p}$&$C_6(q^3-1,Q_{p'})$&$6$\\
\hline
$\chi_{4,0}(k)$&$q^2+q+1$&$\frac{Q_p}{\gcd(Q,k)_p}$&$C_5(Q_{p'})$&$2$\\
\hline
$\chi_{4,1}(k)$&$q^2+q+1$&$\frac{Q_p}{\gcd(Q,k)_p}$&$C_5(Q_{p'})$&$2$\\
\hline
$\chi_{4,2}(k)$&$q^2+q+1$&$\frac{Q_p}{\gcd(Q,k)_p}$&$C_5(Q_{p'})$&$2$\\
\hline
$\chi_{5,0}(k)$&$q-1$&$\frac{Q_p}{\gcd(Q,k)_p}$&$C_3(Q_{p'})$&$2$\\
\hline
$\chi_{5,1}(k)$&$q-1$&$\frac{Q_p}{\gcd(Q,k)_p}$&$C_3(Q_{p'})$&$2$\\
\hline
$\chi_{6,0}(k,l)$&$q-1$&$\frac{Q_p}{\gcd(Q,\lcm(k,l))_p}$&$C_3(Q_{p'})$&$\{4,6\}$\\
\hline
$\chi_{7,0}(k)$&$q+1$&$\frac{Q_p}{\gcd(Q,k)_p}$&$C_{15}(q^3+1,Q_{p'})$&$6$\\
\hline
$\chi_{7,1}(k)$&$q+1$&$\frac{Q_p}{\gcd(Q,k)_p}$&$C_{15}(q^3+1,Q_{p'})$&$6$\\
\hline
$\chi_{9,0}(k)$&$q^2-q+1$&$\frac{Q_p}{\gcd(Q,k)_p}$&$C_{10}(Q_{p'})$&$2$\\
\hline
$\chi_{9,1}(k)$&$q^2-q+1$&$\frac{Q_p}{\gcd(Q,k)_p}$&$C_{10}(Q_{p'})$&$2$\\
\hline
$\chi_{9,2}(k)$&$q^2-q+1$&$\frac{Q_p}{\gcd(Q,k)_p}$&$C_{10}(Q_{p'})$&$2$\\
\hline
$\chi_{10,0}(k)$&$q+1$&$\frac{Q_p}{\gcd(Q,k)_p}$&$C_7(Q_{p'})$&$2$\\
\hline
$\chi_{10,1}(k)$&$q+1$&$\frac{Q_p}{\gcd(Q,k)_p}$&$C_7(Q_{p'})$&$2$\\
\hline
$\chi_{12,0}(k,l)$&$q^2+q+1$&$\frac{Q_p}{\gcd(Q,\lcm(k,l))_p}$&$C_4(Q_{p'})$&$\{6,8\}$\\
\hline
$\chi_{13,0}(k,l)$&$q^2-q+1$&$\frac{Q_p}{\gcd(Q,\lcm(k,l))_p}$&$C_9(Q_{p'})$&$\{6,8\}$\\
\hline
$\chi_{14,0}(k)$&$q^4-q^2+1$&$\frac{Q_p}{\gcd(Q,k)_p}$&$C_{14}(Q_{p'})$&$4$\\
\hline
$\chi_{15,0}(k,l)$&$q+1$&$\frac{Q_p}{\gcd(Q,\lcm(k,l))_p}$&$C_7(Q_{p'})$&$\{4,6\}$\\
\hline
\end{tabular}\\
$\ $\\
$\ $\\
\begin{tabular}{|c|c|c|c|c|}
\multicolumn{5}{c}{Table 12 \ $p$-elements for ${^2}B_2(q),\ q=2^{2n+1}$}\\
\multicolumn{5}{c}{(for the character table in \cite{CHEVIE} which was computed in \cite{Suz})}\\
\hline
character $\chi$&polynomial $Q$&$p$-power $f'_p$&class of $g$&$\Tilde{\lambda}(\chi(g),f'_p)$\\
\hline
\hline
$\chi_{2,0}(k)$&$q-1$&$\frac{Q_p}{\gcd(Q,k)_p}$&$C_2(Q_{p'})$&$2$\\
\hline
$\chi_{3,0}(k)$&$q+2^{n+1}+1$&$\frac{Q_p}{\gcd(Q,k)_p}$&$C_3(Q_{p'})$&$4$\\
\hline
$\chi_{4,0}(k)$&$q-2^{n+1}+1$&$\frac{Q_p}{\gcd(Q,k)_p}$&$C_4(Q_{p'})$&$4$\\
\hline
\end{tabular}\\
$\ $\\
$\ $\\
\begin{tabular}{|c|c|c|c|c|}
\multicolumn{5}{c}{Table 13 \ $p$-elements for ${^2}F_4(q),\ q=2^{2n+1}$}\\
\multicolumn{5}{c}{(for the character table in \cite{CHEVIE})}\\
\hline
character $\chi$&polynomial $Q$&$p$-power $f'_p$&class of $g$&$\Tilde{\lambda}(\chi(g),f'_p)$\\
\hline
\hline
$\chi_{2,0}(k)$&$q-1$&$\frac{Q_p}{\gcd(Q,k)_p}$&$C_3(Q_{p'})$&$\{2,4\}$\\
\hline
$\chi_{2,3}(k)$&$q-1$&$\frac{Q_p}{\gcd(Q,k)_p}$&$C_3(Q_{p'})$&$\{2,4\}$\\
\hline
$\chi_{3,0}(k)$&$q-1$&$\frac{Q_p}{\gcd(Q,k)_p}$&$C_2(Q_{p'})$&$\{2,4\}$\\
\hline
$\chi_{3,1}(k)$&$q-1$&$\frac{Q_p}{\gcd(Q,k)_p}$&$C_2(Q_{p'})$&$\{2,4\}$\\
\hline
$\chi_{4,0}(k,l)$&$q-1$&$\frac{Q_p}{\gcd(Q,\lcm(k,l))_p}$&$C_2(Q_{p'})$&$\{6,8\}$\\
\hline
$\chi_{6,0}(k)$&$q+1$&$\frac{Q_p}{\gcd(Q,k)_p}$&$C_7(Q_{p'})$&$2$\\
\hline
$\chi_{6,1}(k)$&$q+1$&$\frac{Q_p}{\gcd(Q,k)_p}$&$C_7(Q_{p'})$&$2$\\
\hline
$\chi_{15,0}(k,l)$&$q+1$&$\frac{Q_p}{\gcd(Q,\lcm(k,l))_p}$&$C_6(Q_{p'})$&$A$\\
\hline
$\chi_{16,0}(k)$&$q^2-q+1$&$\frac{Q_p}{\gcd(Q,k)_p}$&$C_{16}(Q_{p'})$&$6$\\
\hline
$\chi_{17,0}(k)$&$Q_1$&$\frac{Q_p}{\gcd(Q,k)_p}$&$C_{17}(Q_{p'})$&$\{10,12\}$\\
\hline
$\chi_{18,0}(k)$&$Q_2$&$\frac{Q_p}{\gcd(Q,k)_p}$&$C_{18}(Q_{p'})$&$\{10,12\}$\\
\hline
\multicolumn{5}{l}{$A=\{2,4,6,8,10,12,14,16,18,20,22,24\}$}\\
\multicolumn{5}{l}{$Q_1=q^2-2^{n+1}q+q-2^{n+1}+1$}\\
\multicolumn{5}{l}{$Q_2=q^2+2^{n+1}q+q+2^{n+1}+1$}\\
\end{tabular}\\
$\ $\\
$\ $\\
\begin{tabular}{|c|c|c|c|c|}
\multicolumn{5}{c}{Table 14 \ $p$-elements for ${^2}G_2(q),\ q=3^{2n+1}$}\\
\multicolumn{5}{c}{(for the character table in \cite{CHEVIE} which was computed in \cite{War})}\\
\hline
character $\chi$&polynomial $Q$&$p$-power $f'_p$&class of $g$&$\Tilde{\lambda}(\chi(g),f'_p)$\\
\hline
\hline
$\chi_{3,0}(k)$&$q-1$&$\frac{Q_p}{\gcd(Q,k)_p}$&$C_3(Q_{p'})$&$2$\\
\hline
$\chi_{4,0}(k)$&$q-3^{n+1}+1$&$\frac{Q_p}{\gcd(Q,k)_p}$&$C_4(Q_{p'})$&$6$\\
\hline
$\chi_{5,0}(k)$&$q+1$&$\frac{Q_p}{\gcd(Q,k)_p}$&$C_5(Q_{p'},2)$&$\{2,3,4,5,6\}$\\
\hline
$\chi_{6,0}(k)$&$q+3^{n+1}+1$&$\frac{Q_p}{\gcd(Q,k)_p}$&$C_6(Q_{p'})$&$6$\\
\hline
\end{tabular}\\

\printbibliography

\end{document}